\documentclass{amsart}

\usepackage[utf8]{inputenc}
\usepackage[T1]{fontenc}

\usepackage{geometry}
\usepackage{ragged2e}

\usepackage[pagebackref, hidelinks]{hyperref}
\hypersetup{
colorlinks=true, 
citecolor = black, 
linkcolor = blue, 
urlcolor = blue, 
anchorcolor = blue
}

\usepackage[all]{xy}
\usepackage{mathtools} 
\usepackage{caption}

\usepackage{tikz, pgf, pgfplots}
\pgfplotsset{compat=newest}

\usetikzlibrary[patterns]
\usetikzlibrary{arrows, arrows.meta}
\usetikzlibrary{calc} 
\def\centerarc[#1](#2)(#3:#4:#5)%
    {\draw[#1] ($(#2)+({#5*cos(#3)},{#5*sin(#3)})$) arc (#3:#4:#5);}

\usetikzlibrary{intersections,through,backgrounds} 
\usetikzlibrary{shapes.geometric} 

\definecolor{grey}{rgb}{0.5,0.5,0.5}
\definecolor{forestgreen}{rgb}{0.133,0.545,0.133}
\definecolor{marron}{rgb}{0.6,0.2,0.}
\definecolor{violet}{rgb}{0.5,0.,0.5}
\definecolor{lightpurple}{rgb}{.6,.2,.6}

\usepackage{amsmath, amsthm, amsfonts}
\usepackage{amscd, amssymb, latexsym}

\usepackage{commands_PSL2K}

\usepackage{comment}
\usepackage{color}
\usepackage{url}

\title{Conjugacy classes in $\PSL_2(\Field)$}

\date{\today}
\author{Christopher-Lloyd Simon}

\begin{document}

\maketitle

\begin{abstract}
    We first describe, over a field $\Field$ of characteristic different from $2$, the orbits for the adjoint actions of the Lie groups $\PGL_2(\Field)$ and $\PSL_2(\Field)$ on their Lie algebra $\Sl_2(\Field)$. 
    While the former are well known, the latter lead to the resolution of generalised Pell-Fermat equations which characterise the corresponding orbit.
    The synthetic approach enables to change the base field, and we illustrate this picture over the fields with three and five elements, in relation with the geometry of the tetrahedral and icosahedral groups. 
    While the results may appear familiar, they do not seem to be covered in such generality or detail by the existing literature.
    
    We apply this discussion to partition the set of $\PSL_2(\Z)$-classes of integral binary quadratic forms into groups of $\PSL_2(\Field)$-classes.
    When $\Field = \C$ we obtain the class groups of a given discriminant. Then we provide a complete description of their partition into $\PSL_2(\Q)$-classes in terms of Hilbert symbols, and relate this to the partition into genera.
    The results are classical, but our geometrical approach is of independent interest as it may yield new insights into the geometry of Gauss composition, and unify the picture over function fields.
    %
    
    Finally we provide a geometric interpretation in the modular orbifold $\PSL_2(\Z)\backslash \H$ for when two points or two closed geodesics correspond to $\PSL_2(\Field)$-equivalent quadratic forms, in terms of hyperbolic distances and angles between those modular cycles.
    These geometric quantities are related to linking numbers of modular knots.
    Their distribution properties could be studied using the geometry of the quadratic lattice $(\Sl_2(\Z), \det)$ but such investigations are not pursued here.
    \end{abstract}

\renewcommand{\contentsname}{Plan of the paper}
\setcounter{tocdepth}{1}
\tableofcontents
\subsection*{Acknowlegements}
This article contains the main results obtained in the chapter 1 of my thesis.
I would thus like to thank my thesis advisors Etienne Ghys and Patrick Popescu-Pampu for their guidance and encouragement; 
as well as Francis Bonahon, Louis Funar, Jean-Pierre Otal and Anne Pichon who refereed and carefully read my work.
I am also grateful to Nicolas Bergeron for sharing stimulating discussions, 
%
and to the members of the arithmetics \& dynamics teams at Penn State for inviting me to present my works in their seminars.
%
I owe Marie Dossin for helping me with the figures in tikz.
%
Finally, I thank the editors and referees of the \href{https://mrr.centre-mersenne.org/}{MRR} for their instructive comments.

\newpage

\section*{Introduction}

\subsection*{Adjoint action of \texorpdfstring{$\PSL_2(\Field)$}{PSL(2;K)} on \texorpdfstring{$\Sl_2(\Field)$}{sl(2;K)}}

Let us work over a field $\Field$ of characteristic different from $2$.
The automorphism group $\PGL_2(\Field)$ of the projective line and its largest simple subgroup $\PSL_2(\Field)$ play a fundamental role in various areas of mathematics.
They appear for instance in algebraic geometry when $\Field=\C$ and in hyperbolic geometry when $\Field=\R$ ; in arithmetics when $\Field=\Q$ and in Galois theory when $\Field = \Z/p$.

The first step to understand (the representation theory of) those linear algebraic groups is to describe their conjugacy classes, and more precisely the adjoint actions on their Lie algebra $\Sl_2(\Field)$.
These actions preserve the Killing form, which is a multiple of the non-degenerate quadratic form $\det \colon \Sl_2(\Field) \to \Field$.
%
%
%
It is well known that $\PGL_2(\Field)$ acts transitively on every level set of $\det$. 
After introducing the cross-ratio $\bir(\Aa,\Bb)\in \Field$ of two elements $\Aa,\Bb \in \Sl_2(\Field)$ with non-zero determinant, we will precise this statement.
%

\begin{Proposition}
Let $\Aa,\Bb\in \Sl(\V)$ have determinant $-\delta \ne 0$ and $\bir(\Aa,\Bb) \notin\{1,\infty\}$.

The matrices $M\in \PGL_2(\Field)$ conjugating $\Aa$ to $\Bb$ have a well defined determinant in the quotient $\Field^\times/\Norm_\Field(\Field[\sqrt{\delta}]^\times)$, and its is equal to the class of $\bir(\Aa,\Bb)$.
\end{Proposition}


In contrast, we will index the $\PSL_2(\Field)$-orbits inside $\{\det = -\delta\}$ by the classes $[\chi]$ in the quotient group $\Field^\times / \Norm_\Field(\Field[\sqrt{\delta}])$, and parametrize each orbit $(\delta,\chi)$ by the solutions in $\Field \times \Field$ to the generalised Pell-Fermat equation $x^2-\delta y^2 = \chi$.
In homological terms, the conjugacy problem in $\PSL_2(\Field)$ has obstructions measured by the group $\Field^\times / \Norm_\Field(\Field[\sqrt{\delta}]^\times)$, and when they vanish the conjugacies form a torsor under the group of units $\{\gamma \in \Field[\sqrt{\delta}] \mid \Norm(\gamma)=1 \}$.

\begin{Theorem}
\label{ThmIntro:Conj-a-b-SL(V)}
Let $\Aa,\Bb\in \Sl_2(\Field)$ have determinant $-\delta \ne 0$ and cross-ratio $\bir(\Aa,\Bb)=4\chi \notin \{1,\infty\}$.
The elements $C\in \SL_2(\Field)$ such that $C\Aa C^{-1}= \Bb$ are parametrized by the Pell-Fermat conic:
\begin{equation*}
    (x,y)\in \Field\times \Field 
    \: \colon\: \quad
    x^2-\delta y^2 = \chi
\end{equation*}
\begin{equation*}
    C(x,y) = x(\Id+\Bb\Aa^{-1})+y(\Aa+\Bb)
\end{equation*}

In particular, $\Aa$ and $\Bb$ are conjugate by an element $C(x,y)\in \SL_2(\Field)$ if and only if $\bir(\Aa,\Bb)$ belongs to the subgroup of norms $\Norm_\Field \Field[\sqrt{\delta}]\subset  \Field^\times$ of the quadratic extension, and by an element $C(x,0)\in \SL_2(\Field) \cap \Field[\{\Aa,\Bb\}]$ if and only if $\bir(\Aa,\Bb)$ belongs to the subgroup of squares $(\Field^\times)^2 \subset  \Field^\times$.
\end{Theorem}



The proofs of these statements can be reduced to elementary linear algebra once we thoroughly understand the geometry of the Lie algebra $\Sl_2(\Field)$ inside the quaternion algebra $\Gl_2(\Field)$.
In short, commutativity rhymes with colinearity whereas anti-commutativity rhymes with orthogonality.

We recall this background material in the first two sections 
, and provide an amusing application in the third to prove an analogue of Ptolemy's identity for quadrilaterals inscribed in the isotropic cone $\Cone$ of $(\Sl_2(\Field),\det)$, relying on a natural quadratic desingularization $\psi \colon \Field^2 \to \Cone$.

\subsection*{Classes of binary quadratic forms}

By polarizing a binary quadratic form on $\Field^2$ with respect to the canonical symplectic form $\det$, we obtain an isomorphism:
\begin{equation*}
    Q=lx^2+mxy+ry^2 
    \in \Qfb(\Field^2)
    \qquad 
    \xleftrightarrow[]{Q(v)=\det(v,\Qq v)} \qquad
    \Qq = \tfrac{1}{2}
    \begin{psmallmatrix}
    -m & -2r\\
    2l & m
    \end{psmallmatrix}
    \in \Sl_2(\Field)
\end{equation*}
between the Poisson algebra $(\Qfb(\Field^2),\disc)$ and the Lie algebra $(\Sl_2(\Field),-4\det)$, conjugating the actions of $\PSL_2(\Field)$ by change of variables and by conjugacy.
After interpreting the values of quadratic forms as the scalar products with elements in the isotropic cone $Q(v)=\langle \Qq, \psi(v) \rangle$, and computing that the products $Q_a(u)Q_b(v)=\langle \Aa, \psi(u) \rangle\langle \Bb, \psi(v) \rangle$ for primitive vectors $u,v\in \Field^2$ are equivalent to the cross-ratio $\bir(Q_a,Q_b)=\bir(\Aa,\Bb)$ in $\Field^\times/\Norm(\Field[\sqrt{\Delta}]^\times)$, we will deduce the following.

\begin{Proposition} \label{PropIntro:ClF(D)-group}
The set $\Cl_\Field(\Delta)$ of $\PSL_2(\Field)$-orbits in $\Qfb(\Field^2)$ with non-square discriminant $\Delta$ embeds into the group $\Field^\times /{\Norm_\Field(\Field[\sqrt{\Delta}]^\times})$ of exponent two, by sending the class of the norm $x^2-\tfrac{\Delta}{4}y^2$ of the $\Field$-extension $\Field[\sqrt{\Delta}]$ to the identity, and using the multiplication of values for composition.
\end{Proposition}

The initial motivation was to understand the space $\Qfb(\Z^2)$ of integral binary quadratic forms $Q(x,y)=lx^2+mxy+ry^2$ up to change of variables by $\PSL_2(\Z)$, and the class groups $\Cl(\Delta)$ of primitive classes with non-square discriminant $\Delta = m^2-4lr \in \Z$ introduced by Gauss in \cite{Gauss_disquisitiones_1807}. We refer to \cite{Cassels_Rational-quadratic-forms_1978, Cox_primes-of-the-form_1997} and \cite{Weil_number-theory-history_1984} for the relevant background and history.


For a field $\Field$ of characteristic $\ne 2$, the extension of scalars $\Z\to \Field$ induces a map $\Qfb(\Z^2) \to \Qfb(\Field^2)$, yielding a group morphism $\Cl(\Delta)\to \Cl_\Field(\Delta)$.
We say that $Q_a,Q_b\in \Qfb(\Z^2)$ are $\Field$-equivalent when they are conjugate by $C\in \PSL_2(\Field)$.
When $\Field\supset \Q$ this implies that they have the same discriminant (not just modulo $(\Field^\times)^\times$), and Theorem \ref{ThmIntro:Conj-a-b-SL(V)} \& Proposition \ref{PropIntro:ClF(D)-group} show that their $\Field$-equivalence is measured by $\bir(Q_a,Q_b)\equiv Q_a(1,0)Q_b(1,0)=l_al_b \in \Field^\times/\Norm(\Field[\sqrt{\Delta}]^\times)$. 

%
%

Thus, when $\Field = \C$ this groups the $\PSL_2(\Z)$-classes of $\Qfb(\Z^2)$ according to their discriminant $\Delta$, and as $\Field$ decreases we obtain finer partitions of the class groups $\Cl(\Delta)$ into $\Field$-classes.
In section \ref{Sec:Qfb} we provide a computable characterisation of $\Q$-equivalence in terms of the Hilbert symbols $(\delta,\chi)_p$ at all primes $p\in \Z$, which measures the obstruction to solving the equation $x^2-\delta y^2 = \chi$ in $\Q_p$.
%

We also deduce from Theorem \ref{ThmIntro:Conj-a-b-SL(V)} a relation between the partition of $\Cl(\Delta)$ into $\Q$-classes and its partition into genera, which are given by the cosets $\Cl(\Delta)/\Cl(\Delta)^2$ modulo the subgroup of squares.

\begin{Theorem}
For all non-square discriminant $\Delta\in \Z$, genus equivalence implies $\Q$-equivalence.

For all fundamental discriminants $\Delta$, genus equivalence is also implied by $\Q$-equivalence.
\end{Theorem}


\subsection*{Arithmetic equivalence of singular moduli \& modular geodesics}

We conclude with a geometric interpretation of $\Field$-equivalence, which one may compare with \cite{Penner_geometry-gauss_1996}.
The modular group $\PSL_2(\Z)$ acts on the upper-half plane $\H\P= \{z\in \C \mid \Im(z)>0\}$ by linear fractional transformations, and the quotient is the modular orbifold $\M=\PSL_2(\Z) \backslash \H\P$. 

Consider primitive integral binary quadratic forms $Q_a,Q_b$ with non-square discriminant $\Delta$,  
and denote $(\alpha',\alpha)$, $(\beta',\beta)$ their roots (which one may order up to simultaneous inversion).

If $\Delta>0$, then $Q_a$ and $Q_b$ are uniquely determined by $\alpha,\beta \in \H\P$, and their $\PSL_2(\Z)$-classes correspond to points $[\alpha],[\beta]\in \M$, called singular moduli. 

\begin{Corollary}
Two complex irrationals $\alpha,\beta \in \Q(\sqrt{\Delta})$ are $\Field$-equivalent if and only if there exists a hyperbolic geodesic arc in $\M$ from $[\alpha]$ to $[\beta]$ whose length $\lambda$ is of the form $\left(\cosh\tfrac{\lambda}{2}\right)^2 = \frac{1}{(2x)^2-\Delta y^2}$ for $x,y\in \Field$, in which case all geodesic arcs from $[\alpha]$ to $[\beta]$ have this property.
\end{Corollary}

If $\Delta < 0$ then $Q_a$ and $Q_b$ correspond to the oriented geodesics $(\alpha',\alpha), (\beta',\beta)$ in $\H\P$ and their $\PSL_2(\Z)$-classes correspond to primitive closed oriented geodesics in $\M$, called modular geodesics.

\begin{Corollary}
Two modular geodesics of the same length $2\sinh^{-1}(\sqrt{\Delta}/2)$ are $\Field$-equivalent if and only if one of the following equivalent conditions hold:

\begin{enumerate}
    \item[$\theta$] There exists one intersection point with angle $\theta\in \,]0,\pi[$ such that
    $\left(\cos \tfrac{\theta}{2}\right)^2 = \frac{1}{(2x)^2-\Delta y^2}$ for $x,y\in \Field$,
    in which case all intersections have this property.
    \item[$\lambda$] There exists one co-oriented ortho-geodesic of length $\lambda$ such that $\left(\cosh \frac{\lambda}{2}\right)^2 = \tfrac{1}{(2x)^2-\Delta y^2}$ for $x,y\in \Field$,
    in which case all such ortho-geodesics have this property.
\end{enumerate}
\end{Corollary}
In conclusion, the $\Q$-equivalence is measured by the geometric quantities $(\cos \tfrac{\theta}{2})^2$ or $\left(\cosh \tfrac{\lambda}{2}\right)^2$ as elements in $\Q^\times \bmod \Norm_\Q \Q(\sqrt{\Delta})$, and their multiplication implies a geometric interpretation for the multiplication of genera.

\newpage




\section{Geometric algebra of \texorpdfstring{$\Gl_2(\Field)$}{gl(2;K)}}
\label{Sec:gl(2)}

Consider a field $\Field$ of characteristic different from $2$ and a $\Field$-vector space $\V$ of dimension $2$.

\subsection*{Involutive algebra}

The $\Field$-algebra $\Gl(\V)$ of linear endomorphisms of $\V$ is isomorphic to $\V\otimes \V^*$ with product defined by $(u\otimes \mu)\cdot(v\otimes \nu)=u\otimes \mu(v)\nu$.
%
%
It is endowed with the canonical linear form $\Tr \colon \Gl(\V) \to \Field$ defined by $\Tr(u\otimes \mu) = \mu(u)$.
%
Let us find a canonical \emph{involution} $M\mapsto M^\#$ on $\Gl(\V)$, that is an anti-commutative linear endomorphism of order two, to deduce a canonical  non-degenerate bilinear form $(M,N)\mapsto \Tr(MN^\#)$ isomorphic to the pairing $\Gl(\V) \times \Gl(\V)^* \to \Field$.

A non-degenerate bilinear form $\omega \colon \V\times \V \to \Field$ is equivalent to an isomorphism $\omega^* \colon \V \to \V^*$. Its associated adjoint involution $\#\colon \Gl(\V)\to \Gl(\V)$ is the composition of $(\omega^*)\otimes (\omega^*)^{-1} \colon \V\otimes \V^* \to \V^* \otimes \V$ with the canonical map switching factors, thus $\# \colon u\otimes \mu \mapsto (\omega^*)^{-1}(\mu) \otimes \omega^*(u)$.
Observe that $\#$ only depends on $\omega$ up to scaling.
The plane $\V$ admits a unique non-degenerate anti-symmetric bilinear form $\omega$ up to scaling as $\Lambda^2\V^*\simeq \Field$, and this defines our canonical involution $\#$.

%
Only after choosing a basis of $\V$ do we have the identifications $\V = \Field^2$ and $\Gl(\V)=\Gl_2(\Field)$.
Then $\omega(u,v) = \det(u,v)$ and the associated adjoint involution on $\Gl_2(\Field)$ is the transpose-comatrix:
\begin{equation*}
    M = \begin{psmallmatrix} a & b \\ c & d \end{psmallmatrix}
    \longmapsto 
    M^\# = \begin{psmallmatrix} d & -b \\ -c & a \end{psmallmatrix}.
\end{equation*}
The fixed subalgebra of $M\mapsto M^\#$ is reduced to the center $\Field \Id$ of $\Gl(\V)$.
Composing $(M,M^\#)$ with addition or multiplication yields the central elements:
\begin{equation*}
    \Tr(M)\Id:= M+M^\#
    \qquad \mathrm{and} \qquad
    \det(M)\Id:=M\times M^\#
\end{equation*}
which recovers the linear \emph{trace} map $\Tr \colon \Gl(\V) \to \Field$, and defines the multiplicative \emph{determinant} map $\det \colon \Gl(\V) \to \Field$.
The involution $\#$ preserves the group $\GL(\V)$ of invertible elements, which consists in those $A \in \Gl(\V)$ such that $\det(A)\in \Field^\times$, in which case $A^{-1} = \det(A)^{-1} A^\#$.

For $A\in \GL(\V)$ and $M\in \Gl(\V)$ we have $(AMA^{-1})^\#=AM^\#A^{-1}$, so the left adjoint linear action of $\GL(\V)$ on $\Gl(\V)$ preserves the involution, whence all the structures which will follow. 

The kernel $\SL(\V)$ of the determinant morphism $\det \colon \GL(\V) \to \Field^\times$ is called the subgroup of units, thus $A\in \SL(\V)\iff \det(A)=1 \iff A^\# = A^{-1}$.
The kernel $\Sl(\V)$ of the trace form is the anti-symmetric part for the involution, thus $a\in \Sl(\V) \iff \Tr(a) = 0 \iff a^\# = -a$.

\subsection*{Quadratic space}

On the vector space $\Gl(\V)$ the determinant is a non degenerate quadratic form, and as $\det(M+N)\Id= (M+N)(M+N)^\#=\left(\det(M)+\Tr(MN^\#)+\det(N)\right)\Id$, its polar symmetric bilinear form is:
\begin{equation*}
    \langle M,N \rangle = \tr(MN^\#) \quad \mathrm{where} \quad \tr(P):=\tfrac{1}{2}\Tr(P)
\end{equation*}

The involution $\#$ has eigenvalues $\pm 1$ and its eigenspaces provide a decomposition \[\Gl(\V) = \Field \Id \oplus \Sl(\V)\] which is orthogonal with respect to the determinant form. 
Thus every element $M\in \Gl(\V)$ splits as the sum of its symmetric and anti-symmetric parts with respect to the involution: 
\begin{equation*}
M= \tr(M)\Id+\pr(M)
\qquad \mathrm{where} \quad
\tr(M)\Id=\tfrac{M+M^\#}{2}
\quad \mathrm{and} \quad
\pr(M):=\tfrac{M-M^\#}{2}.
\end{equation*}
In particular $\det(M)\Id = \tr(M)^2\Id - \pr(M)^2$ which we may write $\det = \tr^2-\pr^2$.

The $4$-dimensional space $\Gl(\V)$, which contains the isotropic cone $\Gl(\V)\setminus \GL(\V)$ defined by $\det(M)=\langle M,M \rangle = 0$, decomposes as the direct sum of the anistropic line $\Field\Id$ and its orthogonal hyperplane $\Sl(\V)$ defined by $\tr(M)=\langle \Id,M \rangle = 0$.
Denote by $\Cone$ the isotropic cone for the determinant restricted to the kernel $\Sl(\V)$ of the trace, in formulae: 
\begin{equation*}
    \Cone=\{M\in \Gl(\V)\mid \langle \Id,M \rangle = 0 = \langle M,M \rangle\} = \{\Aa\in \Sl(\V)\mid \det(\Aa)=0\} 
\end{equation*}


\subsection*{Discriminant}

The relation $M^2-(M+M^\#)M+(MM^\#) = 0$ yields the Cayley-Hamilton identity $\chi_M(M)=0$ for $X^2-\Tr(M)X+\det(M) \in \Field[X]$ the \emph{characteristic polynomial} of $M$.
Hence a non-central element $M\in \Gl_2(\Field)\setminus \Field \Id$ generates a commutative subalgebra $\Field[M]=\Span(\Id,M)$ of dimension $2$, which is isomorphic to the quadratic extension $\Field[X]/(\chi_M)$ of $\Field$ with Galois involution given by the restriction of the involution $\#$.
%

The discriminant of $M\in \Gl(\V)$ is defined as that of its characteristic polynomial, equal to 
\begin{equation*}
    \disc(M)=\Tr(M)^2-4\det(M).
\end{equation*}
In particular $\disc(A)=\Tr(A)^2-4$ for $A\in \SL(\V)$ and $\disc(\Aa)=-4\det(\Aa)$ for $\Aa\in \Sl(\V)$.

We call $M\in \Gl(\V)$ \emph{semi-simple} when $\disc(M)\ne 0$, that is when $\chi_M$ has simple roots in its splitting field.
If these roots belong to $\Field$ then $\Field[M]$ is isomorphic to the direct product $\Field \times \Field$, otherwise $\Field[M]$ is a simple $\Field$-algebra (no proper ideals).
In both cases $\Field[M]$ is a semi-simple $\Field$-algebra (a product of simple algebras).
When $\disc(M)=0$ we have $\chi_M(X)=(X-\lambda)^2$ for $\lambda \in \Field$ so the algebra $\Field[M]$ is not integral (it has zero divisors) as $M-\lambda \Id$ is nilpotent.
%

The discriminant is preserved under the projection $\disc(M) = \disc(\pr M)$, so an element is semi-simple if and only if its projection in $\Sl(\V)$ lies outside the cone $\Cone$.

\subsection*{Projectivization}

In the projective $3$-space $\P(\Gl(\V))$, the point $\P(\Field \Id)$ and the plane $\P(\Sl(\V))$ are mutually polar with respect the non-degenerate quadric $\P(\Gl(\V))\setminus \PGL(\V)$. The point lies off the quadric and its polar plane intersects the quadric transversely along the non-degenerate conic $\P(\Cone)$.
%
%
Geometrically, the conic $\P(\Cone)$ consists in the set of tangency points between the quadric $\P(\{\det = 0\})$ and the pencil of lines through $\P(\Id)$.
%


\begin{figure}[h]
    \centering
    \includegraphics[ width=0.45\textwidth]{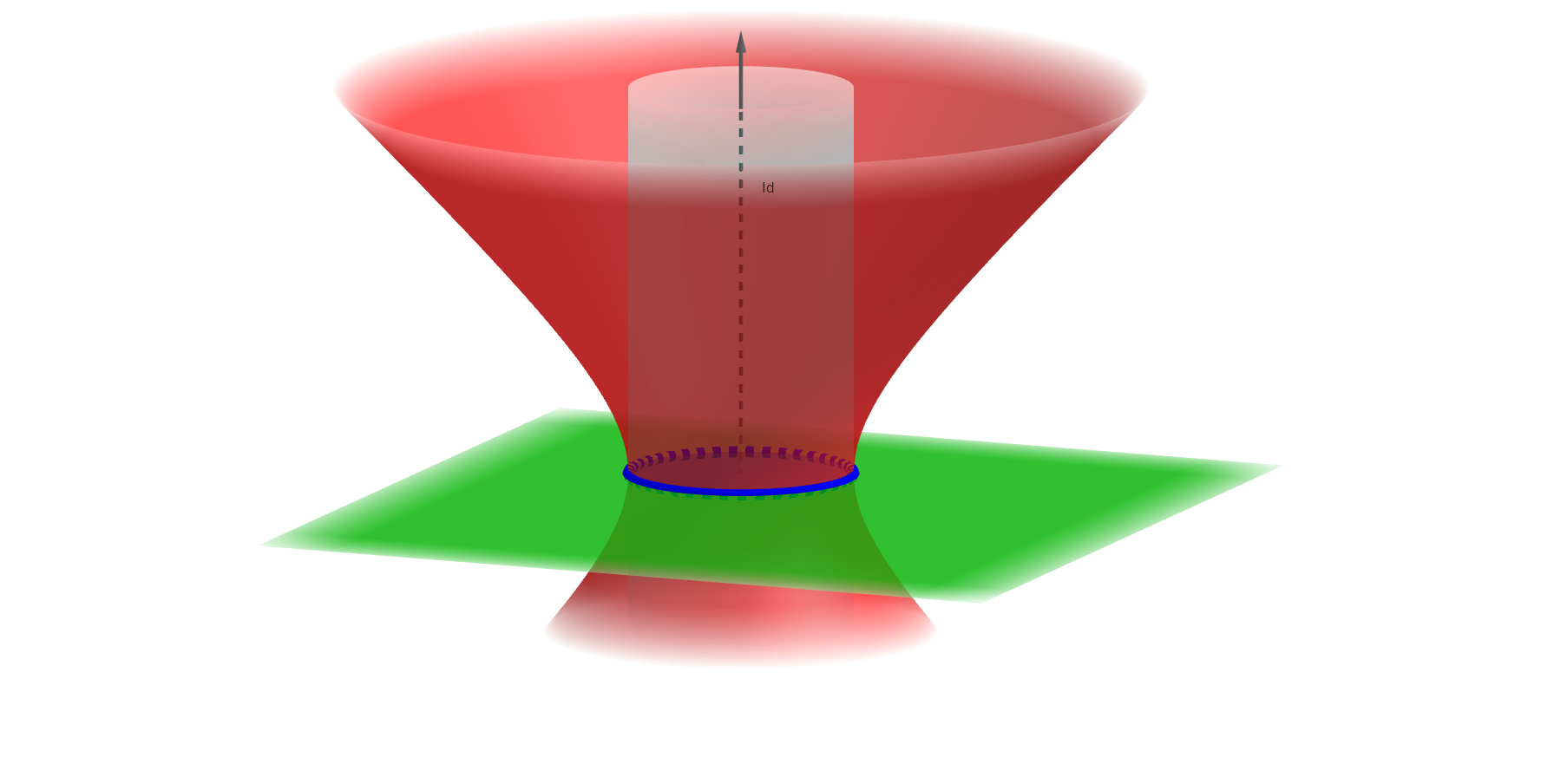}
    \vspace{-0.4cm}
    \caption*{The quadric \textcolor{red}{$\P(\{\det=0\})$} in $\P(\GL(\V))$. The point $\P(\Id)$ lies off the quadric and its polar plane \textcolor{green!64!black}{$\P(\Sl(\V))$} intersects the quadric in the conic \textcolor{blue}{$\P(\Cone)$}.}
    \label{fig:P(Gl(V)-projective-quadric}
\end{figure}

Over $\Field$, the isomorphism types of the quadric $\P(\{\det=0\})$ and of the conic $\P(\Cone)$ are given, in terms of the classes in $\Field^\times/(\Field^\times)^2$ of the diagonal elements appearing in the diagonalisation of the quadratic form $\det$, by $\{1,-1,1,-1\}$ and $\{-1,1,-1\}$.

\begin{Lemma}[Equivariant ruling of the quadric]
\label{Lem:Param-Quadric-Two-Lines}
The map $\Psi \colon M\mapsto (\ker M,\im M)$ defines a bijective algebraic correspondence between the projective quadric $\P(\{\det = 0\})$ and $\P(\V)\times \P(\V)$ sending the projective conic $\P(\Cone)$ to the diagonal $\P(\V)$.
The map $\Psi$ conjugates the adjoint action of $\PGL(\V)$ restricted to $\P(\{\det=0\})$ with its tautological diagonal action on $\P(\V)\times \P(\V)$.
\end{Lemma}




Recall that the action of $\PGL(\V)$ on $\P(\V)$ is simply-transitive on triples of distinct lines. For a symplectic form $\omega$ on $\V$, we define the Maslov index of such a triple $x_1,x_2,x_3\in \P(\V)$ as the element $\omega(\vec{x_1},\vec{x_{2}})\in \Field^\times/(\Field^\times)^2$ where the three vectors $\Vec{x_i}\in x_i \subset \V$ sum up to zero. The level sets of the Maslov index do not depend on the choice of $\omega$. One may show \cite[Proposition 1.39]{CLS_phdthesis_2022} that the action of $\PSL(\V)$ on $\P(\V)$ is simply-transitive on triples of distinct lines with a given Maslov index.


\section{The Lie algebra \texorpdfstring{$\Sl_2(\Field)$}{sl(2;K)}}
\label{Sec:sl(2)}

\subsection*{Orthogonality \& colinearity}

The associative algebra $\Gl(\V)$ inherits the structure of a Lie algebra by taking (half of) the commutator:
\begin{equation*}
    \{M,N\}=\tfrac{1}{2}(MN-NM) 
\end{equation*}
and as $\{M,N\}=\{\pr M,\pr N\}$ we may quotient by its center $\Field \Id$ to find a Lie bracket on $\Sl(\V)$.
From a geometric perspective, the kernel $\Sl(\V)$ of the trace form is a $3$-dimensional Lie algebra, that of the Lie group $\SL(\V)$, kernel of the determinant morphism.

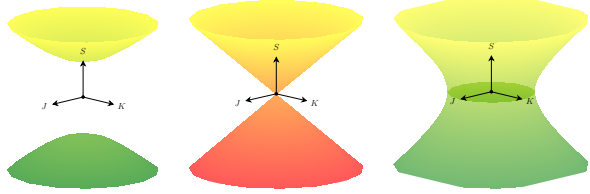
\begin{figure}[h]
    \centering
    \scalebox{0.8}{
\pgfplotsset{compat = newest}

\scalebox{0.2}{
\begin{tikzpicture}
\begin{axis}[
    hide axis,
    view={135}{16},
	xlabel = {$x$},
	ylabel = {$y$},
	zlabel = {$z$},
	legend cell align={left},
	ymin = -2, ymax = 2,
	xmin = -2, xmax = 2,
	xtick={-5,-4,...,5},
	ytick={-5,-4,...,5},
	ztick={-5,-4,...,5},
	scale = 3,
	zmin = -2, zmax = 2,
	z buffer = sort,
]

\addplot3[
    surf,
    shader = interp,
    samples = 16,
    samples y = 17,
    domain  = 0:2*pi,
    domain y = 0:2,
    variable = \t,
    opacity = 0.65,
    colormap/greenyellow,
]
(
    {cos(deg(\t))*sqrt(y^2-1)},
    {(sin(deg(\t))*sqrt(y^2-1)},
    {y}
);

\addplot3[
    surf,
    shader = interp,
    samples = 16,
    samples y = 17,
    domain  = 0:2*pi,
    domain y = 0:2,
    variable = \t,
    opacity = 0.65,
    colormap/greenyellow,
]
(
    {cos(deg(\t))*sqrt(y^2-1)},
    {(sin(deg(\t))*sqrt(y^2-1)},
    {-y}
);

\draw [line width=1.5pt,color=black, fill=black] (0,0,0) circle (3pt);
\draw [line width=2pt, color=black, -{Stealth[width=5mm,length=5mm]}] (0,0,0) -- (1,0,0);
\draw[color=black] (1.25,0,0) node {\Huge$J$};
\draw [line width=2pt, color=black, -{Stealth[width=5mm,length=5mm]}] (0,0,0) -- (0,1,0);
\draw[color=black] (0,1.25,0) node {\Huge$K$};
\draw [line width=2pt, color=black, -{Stealth[width=5mm,length=5mm]}] (0,0,0) -- (0,0,1);
\draw[color=black] (0,0,1.25) node {\Huge$S$};

\end{axis}
\end{tikzpicture}
}
\quad
\scalebox{0.2}{
\begin{tikzpicture}
\begin{axis}[
    hide axis,
    view={135}{15},
	xlabel = {$x$},
	ylabel = {$y$},
	zlabel = {$z$},
	legend cell align={left},
	ymin = -2, ymax = 2,
	xmin = -2, xmax = 2,
	xtick={-5,-4,...,5},
	ztick={-6,-5,...,6},
	scale = 3,
	zmin = -2, zmax = 2,
	z buffer = sort,
]
\addplot3[
    surf,
    shader = interp,
    samples = 16,
    samples y = 17,
    domain  = 0:2*pi,
    domain y = 0:2,
    variable = \t,
    opacity = 0.65,
    colormap/redyellow,
]
(
    {cos(deg(\t))*y},
    {(sin(deg(\t))*y},
    {y}
);

\addplot3[
    surf,
    shader = interp,
    samples = 16,
    samples y = 17,
    domain  = 0:2*pi,
    domain y = 0:2,
    variable = \t,
    opacity = 0.65,
    colormap/redyellow,
]
(
    {cos(deg(\t))*y},
    {(sin(deg(\t))*y},
    {-y}
);

\draw [line width=1.5pt,color=black, fill=black] (0,0,0) circle (3pt);
\draw [line width=2pt, color=black, -{Stealth[width=5mm,length=5mm]}] (0,0,0) -- (1,0,0);
\draw[color=black] (1.25,0,0) node {\Huge$J$};
\draw [line width=2pt, color=black, -{Stealth[width=5mm,length=5mm]}] (0,0,0) -- (0,1,0);
\draw[color=black] (0,1.25,0) node {\Huge$K$};
\draw [line width=2pt, color=black, -{Stealth[width=5mm,length=5mm]}] (0,0,0) -- (0,0,1);
\draw[color=black] (0,0,1.25) node {\Huge$S$};

\end{axis}
\end{tikzpicture}
}
\quad
\scalebox{0.2}{
\begin{tikzpicture}
\begin{axis}[
    hide axis,
    view={135}{16},
	xlabel = {$x$},
	ylabel = {$y$},
	zlabel = {$z$},
	legend cell align={left},
	ymin = -2, ymax = 2,
	xmin = -2, xmax = 2,
	xtick={-5,-4,...,5},
	ztick={-6,-5,...,6},
	scale = 3,
	zmin = -2, zmax = 2,
	z buffer = sort,
]

\addplot3[
    surf,
    shader = interp,
    samples = 16,
    samples y = 17,
    domain  = 0:2*pi,
    domain y = 0:2,
    variable = \t,
    opacity = 0.55,
    colormap/greenyellow,
]
(
    {cos(deg(\t))*sqrt(y^2+1)},
    {(sin(deg(\t))*sqrt(y^2+1)},
    {y}
);

\addplot3[
    surf,
    shader = interp,
    samples = 16,
    samples y = 17,
    domain  = 0:2*pi,
    domain y = 0:2,
    variable = \t,
    opacity = 0.55,
    colormap/greenyellow,
]
(
    {cos(deg(\t))*sqrt(y^2+1)},
    {(sin(deg(\t))*sqrt(y^2+1)},
    {-y}
);

\draw [line width=1.5pt,color=black, fill=black] (0,0,0) circle (3pt);
\draw [line width=2pt, color=black, -{Stealth[width=5mm,length=5mm]}] (0,0,0) -- (1,0,0);
\draw[color=black] (1.25,0,0) node {\Huge$J$};
\draw [line width=2pt, color=black, -{Stealth[width=5mm,length=5mm]}] (0,0,0) -- (0,1,0);
\draw[color=black] (0,1.25,0) node {\Huge$K$};
\draw [line width=2pt, color=black, -{Stealth[width=5mm,length=5mm]}] (0,0,0) -- (0,0,1);
\draw[color=black] (0,0,1.25) node {\Huge$S$};

\end{axis}
\end{tikzpicture}
}
}
    \caption*{The level surfaces $1,0,-1$ of $\det$ in $\Sl_2(\Q)$.}
    \label{fig:levels-det-sl2-basis}
\end{figure}

From now on we will focus on the geometry of $\Sl(\V)$ with the restricted scalar product $\langle\Aa,\Bb\rangle$ and Lie bracket $\{\Aa,\Bb\}$.
For $a,b \in \Sl(\V)$, the decomposition $ab = \tr(ab)\Id + \pr(ab)$ rewrites as 
\begin{equation*}
ab = -\langle a , b\rangle \Id +\{a,b\}
\end{equation*} 
thus \emph{$a,b$ are orthogonal if and only if they anticommute in which case $\{a,b\}=ab$}.
The Jacobi relation implies that $\{a,b\} \perp \Span(a,b)$ for all $a,b\in \Sl(\V)$. (Following \cite{Arnold_Jacobi-Lie_2005, Ivanov_Arnold-Jacobi-Lie_2011}, this yields a geometric interpretation of the Jacobi relation as the orthocenter theorem for triangles in $\P(\Sl(\V))$.)

The Killing form associated to the bracket is proportional to the scalar product: 
\begin{equation*}
    -\tfrac{1}{8}\Tr(c\mapsto 2\{a,2\{b,c\}\})=-\tr(c\mapsto \{a,\{b,c\}\})=-\tr(ab)=\tr(ab^\#)=\langle a,b\rangle
\end{equation*}
The non-degeneracy of the Killing form implies that of the Lie bracket.
It also implies that \emph{$a,b$ are colinear if and only if $\{a,b\}=0$}.
Hence the quantity $[a,b,c]:= \langle \{a,b\},c\rangle = \tfrac{1}{2}(\tr(bac)-\tr(abc))$ defines a volume form on $\Sl(\V)$, that is an alternate non-degenerate trilinear form over $\Field$.

\begin{figure}[h]
    \centering
    \hspace{-2cm}
    \includegraphics[width = 0.42\textwidth]{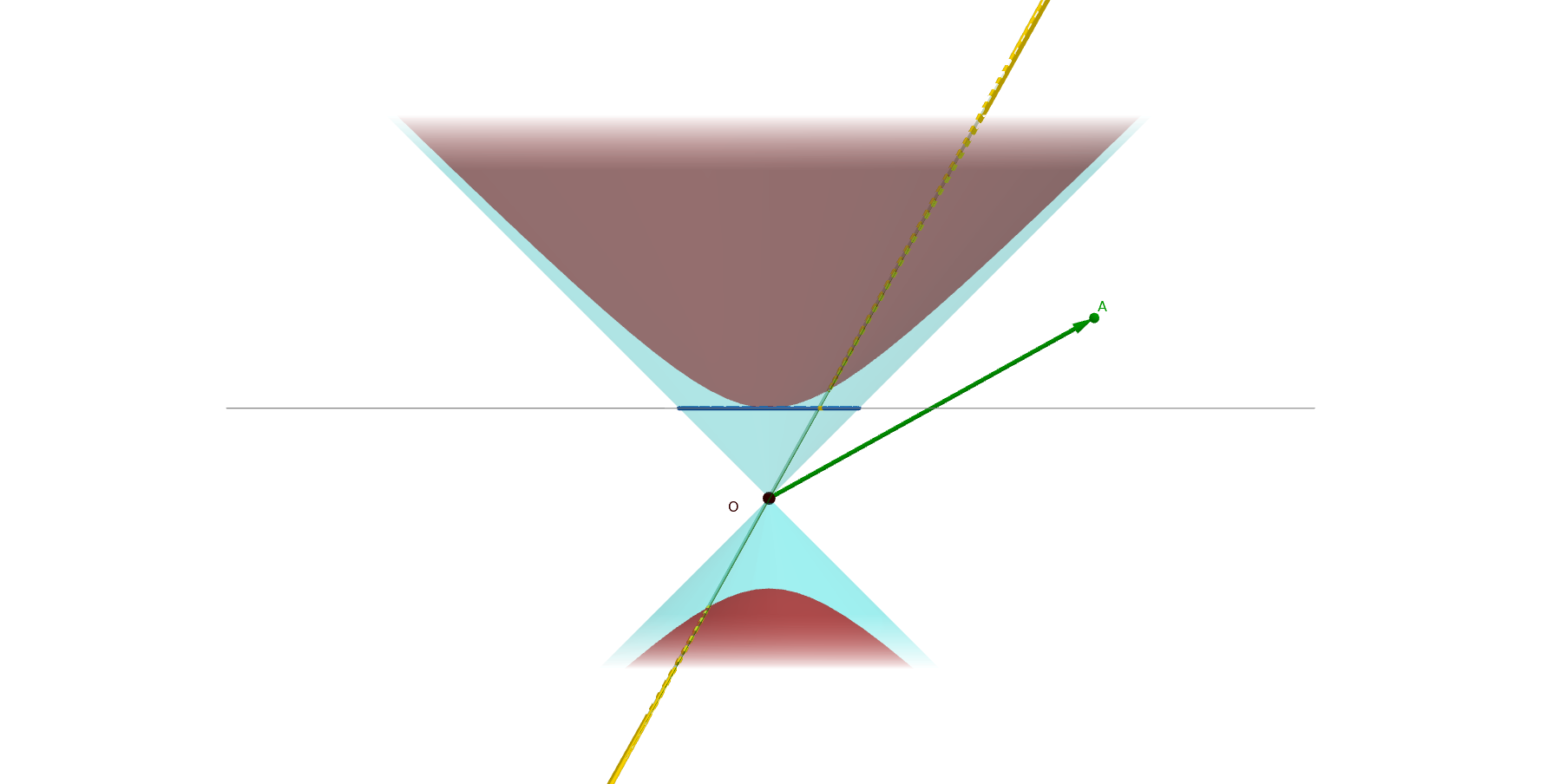}
    \hspace{-2cm}
    \includegraphics[width = 0.42\textwidth]{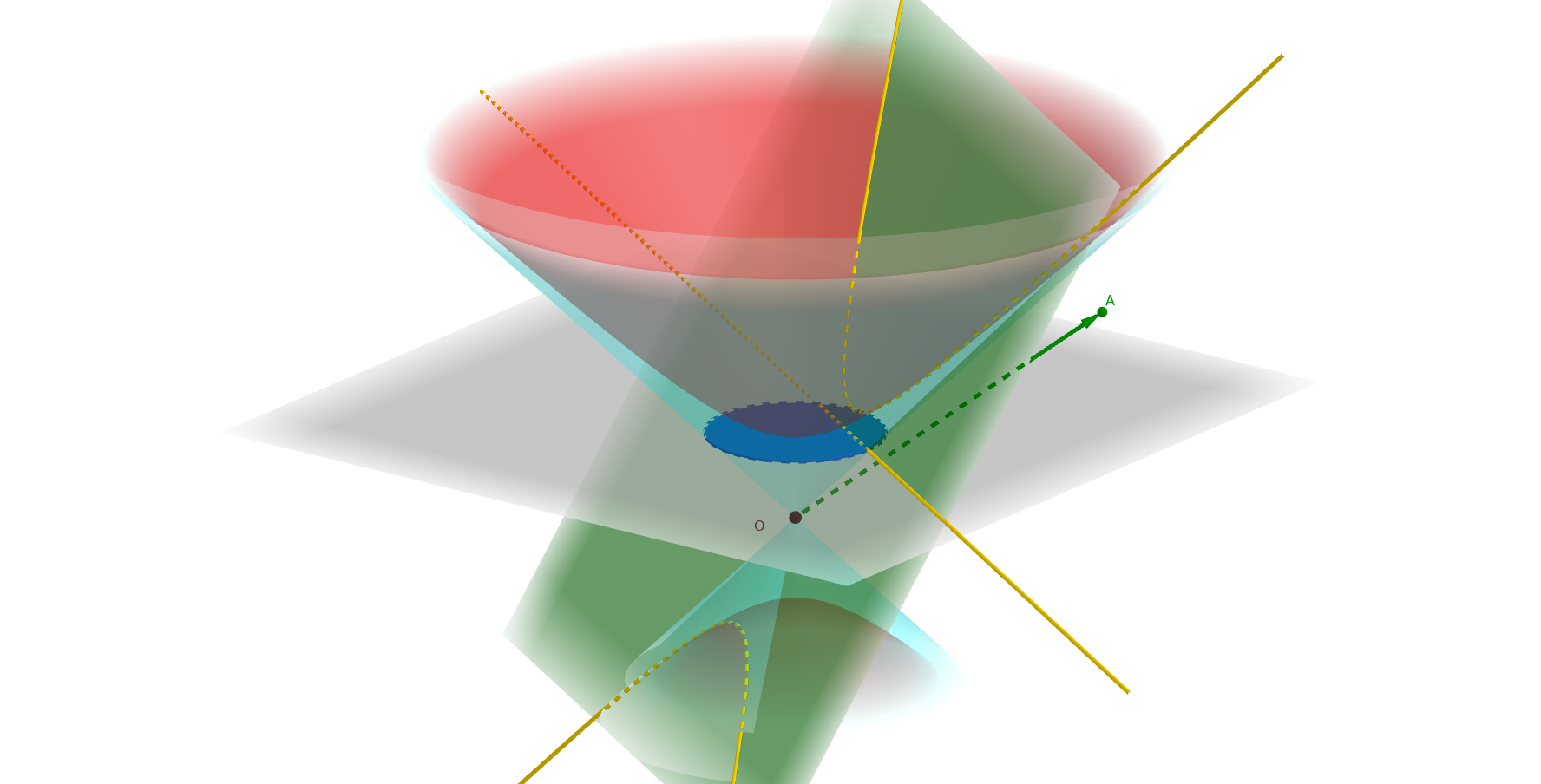}
    \hspace{-2cm}
    \includegraphics[width = 0.42\textwidth]{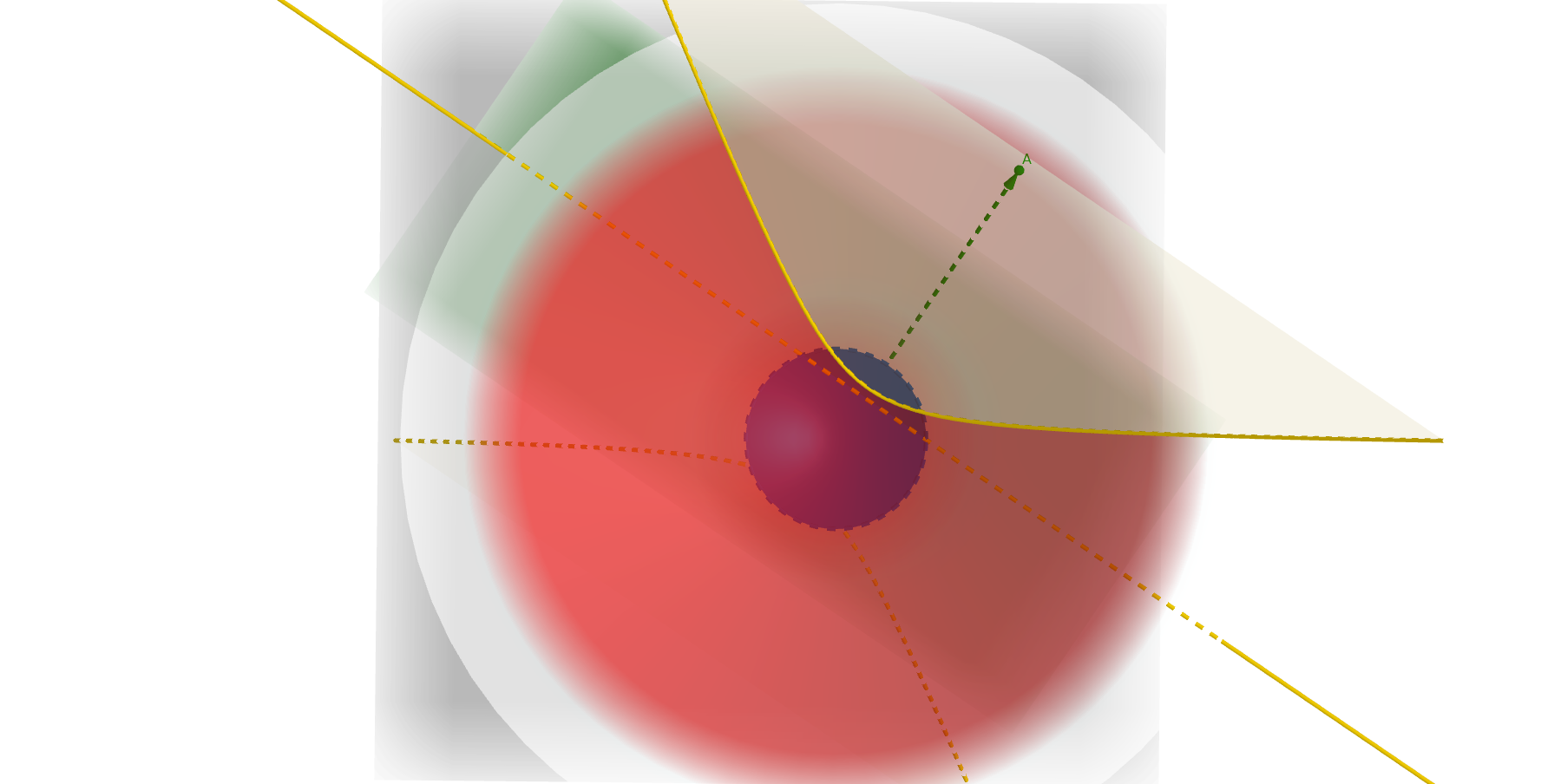}
    \caption*{The projectivization $\P \colon \Sl_2(\Q) \to \Q\P^2$ sends an orthogonal pair $(\Aa,\Aa^\perp)$ to a polar pair (here $\Aa\notin \Cone$ so $a^\perp$ is not tangent to $\Cone$).}
    \label{fig:cayley-klein-projection}
\end{figure}

When $a,b\in \Sl(\V)$ are not colinear, they span the plane $\Span(a,b)=\{a,b\}^\perp$.
In the projective plane $\P(\Sl(\V))$ the line through the distinct points $\P(a),\P(b)$ is polar to the point $\P(\{a,b\})$ with respect to the conic $\P(\Cone)$.
Their relative position with respect to the conic is given by the discriminant $\disc \{a,b\} \in \Field/(\Field^\times)^2$ of the quadratic form $\det$ restricted to the plane $\Span(a,b)$.

\begin{figure}[h]
    \centering
    \scalebox{0.34}{
\begin{tikzpicture}[line cap=round,line join=round,>=triangle 45,x=4.0cm,y=4.0cm]
\clip(-1.1,-1.1) rectangle (1.7,1.1);
\draw [line width=1.2pt] (0,0) circle (4cm);
\draw [line width=2pt,color=violet,domain=-1.1:1.6] plot(\x,{(--1.48-1.8*\x)/-0.4});
\draw [line width=1.2pt,dash pattern=on 2pt off 2pt,color=red,domain=-1.1:1.6] plot(\x,{(-0.74--0.68*\x)/-0.3});
\draw [line width=1.2pt,dash pattern=on 2pt off 2pt,color=blue,domain=-1.1:1.6] plot(\x,{(--0.74-0.49*\x)/-0.56});
\begin{scriptsize}
\fill [color=red] (1,0.8) circle (2.0pt);
\draw[color=red] (1.15,0.85) node {\Huge$a_-$};
\fill [color=blue] (0.6,-1) circle (2.0pt);
\draw[color=blue] (0.77,-0.99) node {\Huge$b_-$};
\fill [color=violet] (1.211,-0.27) circle (2.5pt);
\draw[color=violet] (1.5,-0.27) node {\Huge$\{\textcolor{red}{a},\textcolor{blue}{b}\}$};
\fill [color=red] (0.82,-0.02) circle (2.0pt);
\draw[color=red] (0.67,0.05) node {\Huge$a_+$};
\fill [color=blue] (0.745,-0.35) circle (2.0pt);
\draw[color=blue] (0.6,-0.31) node {\Huge$b_+$};
\fill [color=blue] (0.655,-0.76) circle (2.0pt);
\draw[color=blue] (0.53,-0.7) node {\Huge$b_0$};
\fill [color=red] (0.915,0.41) circle (2.0pt);
\draw[color=red] (0.75,0.43) node {\Huge$a_0$};
\end{scriptsize}
\end{tikzpicture}
}
    \hfill
    \scalebox{0.34}{
\begin{tikzpicture}[line cap=round,line join=round,>=triangle 45,x=4.0cm,y=4.0cm]
\clip(-1.2,-1.2) rectangle (1.8,1.2);
\draw [line width=1.2pt] (0,0) circle (4cm);
\draw [line width=2pt,color=violet,domain=-1.2:1.8] plot(\x,{(--1.36-1.02*\x)/-0.33});
\draw [line width=1pt,dotted,color=red,domain=-1.2:1.8] plot(\x,{(-1.31--0.34*\x)/-1.27});
\draw [line width=1pt,dotted,color=red,domain=-1.2:1.8] plot(\x,{(--1.31-1.13*\x)/-0.67});
\draw [line width=1pt,dotted,color=blue,domain=-1.2:1.8] plot(\x,{(-0.77--0.73*\x)/-0.26});
\draw [line width=1pt,dotted,color=blue,domain=-1.2:1.8] plot(\x,{(--0.77-0.43*\x)/-0.64});
\draw [line width=1.2pt,dash pattern=on 3pt off 3pt,color=violet,domain=-1.2:1.8] plot(\x,{(--0.96-1.47*\x)/0.6});
\draw [line width=1.2pt,dash pattern=on 3pt off 3pt,color=violet,domain=-1.2:1.8] plot(\x,{(--0.97-1.16*\x)/-0.39});
\begin{scriptsize}
\fill [color=red] (1.53,0.62) circle (2.0pt);
\draw[color=red] (1.65,0.46) node {\Huge$a_-$};
\fill [color=blue] (1.2,-0.4) circle (2.0pt);
\draw[color=blue] (1.42,-0.45) node {\Huge$b_-$};
\fill [color=violet] (0.755,-0.245) circle (2.5pt);
\draw[color=violet] (0.42,-0.245) node {\Huge$\{\textcolor{red}{a},\textcolor{blue}{b}\}$};
\end{scriptsize}
\end{tikzpicture}
}
    \hfill
    \scalebox{0.34}{
\begin{tikzpicture}[line cap=round,line join=round,>=triangle 45,x=4.0cm,y=4.0cm]
\clip(-1.2,-1.1) rectangle (1.5,1.1);
\draw [line width=1.2pt] (0,0) circle (4cm);
\draw [line width=1.6pt,color=violet] (1,-1.1) -- (1,1.1);
\begin{scriptsize}
\fill [color=red] (1,0.8) circle (2.0pt);
\draw[color=red] (1.18,0.8) node {\Huge$a_-$};
\fill [color=blue] (1,-0.6) circle (2.0pt);
\draw[color=blue] (1.18,-0.61) node {\Huge$b_-$};
\fill [color=violet] (1,0) circle (2.5pt);
\draw[color=violet] (1.25,0.) node {\Huge$\{\textcolor{red}{a},\textcolor{blue}{b}\}$};
\end{scriptsize}
\end{tikzpicture}
}
    \caption*{In the projective plane $\P(\Sl_2(\Q))$ with the conic $\P(\Cone)$: various configurations of the line $(\P(a),\P(b))$ and its pole $\P(\{a,b\})$.}
    \label{fig:Line-[a,b]-versus-cone-X}
\end{figure}

\subsection*{Subalgebras \& commutants}

Elements $a,b\in \Sl(\V)$ generate an associative subalgebra $(\Field[a,b],\cdot)$ of $(\Gl(\V),\cdot)$ and generate a Lie subalgebra $(\Lie(a,b),\{\cdot,\cdot\})$ of $(\Sl(\V),\{\cdot,\cdot\})$.
Since the Lie bracket equals half the commutator of the associative product, we clearly have $\Field[a,b]\supset \Field \Id \oplus \Lie(a,b)$.

\begin{Proposition}[Subalgebras]
\label{Prop:Lie[a,b]}

Let $a,b\in \Sl(\V)$.
In terms of the underlying vector spaces we have $\Field[a,b]=\Field \Id \oplus \Lie(a,b)$ and $\Lie(a,b)=\Span(a,b,\{a,b\})$.

There are four possibilities for the isomorphism type of $\Lie(a,b)$ given by the relative position of $a,b,\{a,b\}$ with respect to the isotropic cone $\Cone \subset \Sl(\V)$, as follows.

\begin{itemize}
    \item[0] If $a=0=b$, then $\Lie(a,b)=\{0\}$.

    \item[1] If $\{a,b\}=0$ but $\det(a)\ne 0$ then $\Lie(a,b)$ is the abelian Lie algebra of $\dim = 1$. 

    \item[2] If $\{a,b\}\ne 0$ but $\det\{a,b\}=0$ then $\Lie(a,b)$ is the affine Lie algebra of $\dim = 2$. 

    \item[3] If $\det\{a,b\}\ne0$ then $\Lie(a,b)=\Sl(\V)$.
\end{itemize}


Over every field $\Field \nsupseteq \Z/2$, each of these cases can be realised by choosing $a,b$ appropriately.
\end{Proposition}

The previous Proposition could have been formulated for $M,N\in \Gl(\V)$, since they generate the same associative and Lie algebras as their projections $a, b \in \Sl(\V)$.

\begin{Proposition}[Commutants]
\label{Prop:Stab_GL_sl-X}
Consider the adjoint actions of the groups $\GL(\V)$ and $\SL(\V)$ on the space $\Sl(\V)$ and its projectivisation $\P(\Sl(\V))$.
%
Let $p\in \Sl(\V) \setminus \Cone$.

The stabilizer of $p$ under $\GL(\V)$ is $(\Field[p])^\times$, that is the complement of the degenerate conic $x^2+y^2\det(p)=0$ in the plane $\Field[p]=\{x\Id +yp\}$. 

The stabilizer of $p$ under $\SL(\V)$ is $(\Field[p])^\times\cap \SL(\V)$, that is the conic $x^2+y^2\det(p)=1$ in the plane $\Field[p]=\{x\Id +yp\}$, which is non-degenerate except when $p\in \Cone$.

The stabiliser of $\P(p)$ under $\GL(\V)$ is the $\Z/2$-graded subgroup $(\Field[p])^\times \sqcup (\Field[p]^\perp)^\times$ formed by the union of the complements of two degenerate conics.

The stabiliser of $\P(p)$ under $\SL(\V)$ is the $\Z/2$-graded subgroup $(\Field[p]\cap \SL(\V)) \sqcup (\Field[p]^\perp \cap \SL(\V))$ formed by the union of two non-degenerate conics isomorphic to $x^2+y^2\det(p)=\pm 1$.

(The conic $\Field[p]^\perp \cap \SL(\V)$ is isomorphic to $x^2+y^2\det(p)=-1$ but does not have a privileged parametriztation and may have no $\Field$-points.)
\end{Proposition}


\begin{proof}[Proof]
For $C\in \GL(\V)$, if $CpC^{-1}$ is proportional to $p$, then either it equals $+p$ in which case $\{\pr C, p\}=0$ and $C\in \Field[p]$, or else it equals $-p$ in which case $\langle \pr C, p\rangle =0$ and $\pr C \in \Field[p]^\perp$.

Recall that the quadratic space $(\Gl(\V),\det)$ is isomorphic to a sum of two hyperbolic planes. Hence the restrictions of $\det$ to the summands of the decomposition $\Gl(\V) = \Field[p]\oplus \Field[p]^\perp$ must have opposite Witt classes, and the former is $(1,\det(p))$ so the latter is equivalent to $(-1,-\det(p))$.
\end{proof}

We refer to \cite[1.28]{CLS_phdthesis_2022} for analogous descriptions of the stabilizers of $p\in \Cone\setminus \{0\}$.

\subsection*{Cosine and cross-ratio}

For $\Aa \in \Sl(\V)\setminus \Cone$, choose a square root of $\delta:= -\det(\Aa)$ and extend the scalars to the field $\Field'=\Field[\sqrt{\delta}]$.
%
%
The tautological action of $\Aa$ on the plane $\Field'^2$ has two eigendirections for the eigenvalues $\pm \sqrt{\delta}$.
These lines are mapped by $\Psi \otimes \Field'$ to the intersection of the cone $\Cone\otimes \Field'$ with the orthogonal plane $\Aa^\perp$.
We deduce an ordered pair of points $\alpha',\alpha\in \Field'\P^1$.

We may now define and relate the cosine $\cos(\Aa,\Bb)$ \& cross-ratio $\bir(\Aa,\Bb)$ of $\Aa,\Bb \in \Sl(\V)\setminus \Cone$.
These equivalent quantities, together with the discriminants $\disc(\Aa) \& \disc(\Bb)$, are the only $\PGL(\V)$-invariants for a pair of elements in $\Sl(\V)\setminus \Cone$. 

\begin{Lemma}
\label{Lem:bir-cos}
For $\Aa,\Bb \in \Sl(\V)\setminus \Cone$, if we choose a square root of $\det(\Aa\Bb)$ then we may define their cosine $\cos(\Aa,\Bb) \in \Field[\sqrt{\det(\Aa\Bb)}]$:
\begin{equation*}
    \cos(\Aa,\Bb):= \frac{\langle \Aa,\Bb \rangle}{\sqrt{\langle \Aa,\Aa \rangle \langle \Bb,\Bb \rangle}} 
    = \frac{-\tfrac{1}{2}\Tr(\Aa\Bb)}{\sqrt{\det(\Aa\Bb)}} 
\end{equation*}
and we may order their polar points $\P(\Aa^\perp \cap \Cone)=\{\alpha',\alpha\}$ and $\P(\Bb^\perp \cap \Cone)=\{\beta',\beta\}$ up to simultaneous inversion, so as to define their cross-ratio $\bir(\Aa,\Bb)\in \Field[\sqrt{\det(\Aa\Bb)}]$:
\begin{equation*}
    \bir(\Aa,\Bb):= 
    \bir(\alpha',\alpha;\beta',\beta) =
    \frac{(\alpha-\alpha')(\beta-\beta')}{(\alpha-\beta')(\beta-\alpha')} 
\end{equation*}
For a same choice of $\sqrt{\det(\Aa\Bb)}$, these quantities are related by:
\begin{equation*}
    \frac{1}{\bir(\Aa,\Bb)}
    =\frac{1+\cos(\Aa,\Bb)}{2}
\end{equation*}
\end{Lemma}

\begin{Remark}
When $\det(\Aa)=\det(\Bb)$, this common value yields a canonical choice for $\sqrt{\det(\Aa\Bb)}$.
\end{Remark}

\begin{Remark}
For $\Aa,\Bb \in \Sl(\V)\setminus \Cone$:
$\det\{\Aa,\Bb\}= 0 \iff \cos(\Aa,\Bb)^2 = 1 \iff \bir(\Aa,\Bb)\in\{1,\infty\}$.
\end{Remark}

\begin{figure}[h]
    \centering
    \scalebox{0.32}{
\begin{tikzpicture}[line cap=round,line join=round,>=triangle 45,x=4.0cm,y=4.0cm]
\clip(-1.2,-1.2) rectangle (1.3,1.3);
\draw [shift={(0.41,0.5)},line width=1.8pt,color=red,fill=red,fill opacity=0.1] (0,0) -- (-40.55:0.12) arc (-40.55:67.61:0.12) -- cycle;
\draw [shift={(0.41,0.5)},line width=1.8pt,color=red,fill=red,fill opacity=0.1] (0,0) -- (139.45:0.12) arc (139.45:247.61:0.12) -- cycle;
\draw [line width=1.2pt] (0,0) circle (4cm);
\draw [line width=2pt,color=blue] (-0.15,0.99)-- (1,0);
\draw [line width=2pt,color=blue] (-0.2,-0.98)-- (0.55,0.83);
\draw [line width=1.2pt,dash pattern=on 3pt off 3pt,color=blue] (-0.15,0.99)-- (0.55,0.83);
\draw [line width=1.2pt,dash pattern=on 3pt off 3pt,color=blue] (-0.2,-0.98)-- (1,0);
\begin{scriptsize}
\fill [color=blue] (-0.15,0.99) circle (2.0pt);
\draw[color=blue] (-0.3,1.13) node {\Huge$\alpha'$};
\fill [color=blue] (1,0) circle (2.0pt);
\draw[color=blue] (1.2,-0.03) node {\Huge$\alpha$};
\fill [color=blue] (-0.2,-0.98) circle (2.0pt);
\draw[color=blue] (-0.28,-1.12) node {\Huge$\beta'$};
\fill [color=blue] (0.55,0.83) circle (2.0pt);
\draw[color=blue] (0.78,0.92) node {\Huge$\beta$};
\draw[color=red] (0.68,0.55) node {\Huge$\theta$};
\end{scriptsize}
\end{tikzpicture}
}
    \qquad
    \scalebox{0.32}{
\begin{tikzpicture}[line cap=round,line join=round,>=triangle 45,x=4.0cm,y=4.0cm]
\clip(-1.2,-1.2) rectangle (1.3,1.3);
\draw [line width=1.2pt] (0,0) circle (4cm);
\draw [line width=1.2pt,dash pattern=on 3pt off 3pt,color=blue] (-0.15,0.99)-- (1,0);
\draw [line width=1.2pt,dash pattern=on 3pt off 3pt,color=blue] (-0.2,-0.98)-- (0.55,0.83);
\draw [line width=2pt,color=blue] (-0.15,0.99)-- (0.55,0.83);
\draw [line width=2pt,color=blue] (-0.2,-0.98)-- (1,0);
\draw [line width=2pt,color=red] (0.68,-0.27)-- (0.28,0.89);
\begin{scriptsize}
\fill [color=blue] (-0.15,0.99) circle (1.5pt);
\draw[color=blue] (-0.32,1.08) node {\Huge$\alpha'$};
\fill [color=blue] (1,0) circle (2.0pt);
\draw[color=blue] (1.19,0.05) node {\Huge$\beta$};
\fill [color=blue] (-0.2,-0.98) circle (2.0pt);
\draw[color=blue] (-0.32,-1.1) node {\Huge$\beta'$};
\fill [color=blue] (0.55,0.83) circle (2.0pt);
\draw[color=blue] (0.87,0.84) node {\Huge$\alpha$};
\draw[color=red] (0.62,0.57) node {\Huge$\lambda$};
\end{scriptsize}
\end{tikzpicture}
}
    \caption*{Angle at intersection $\frac{1}{\bir} =\left(\cos \tfrac{\theta}{2}\right)^2$. Length of ortho-geodesic $\frac{1}{\bir} =\left(\cosh \tfrac{\lambda}{2}\right)^2$.}
    \label{fig:bir-cos}
\end{figure}



\section{Ptolemy's theorem for quadrilaterals inscribed in \texorpdfstring{$\P(\Cone)$}{P(X)}}
\label{Sec:Ptolemy}


\subsection*{Parametrizing the cone}

Choose a symplectic form $\omega$ on $\V$ and consider the associated quadratic map $\psi \colon \V \to \Gl(\V)$ given by $\psi(v)= -v\otimes \omega^*(v)=\omega(\cdot,v)v$.
Thus $\psi(v)$ is nilpotent with image $\Field v$.
Consequently, $\psi(u)\psi(v) \in \Gl(\V)$ equals $-\omega(u,v)^2$ times the projection on $\Field u$ parallel to $\Field v$, whence $\{\psi(u),\psi(v)\}\in \Sl(\V)$ equals $-\omega(u,v)^2$ times the symmetry with respect $\Field u$ parallel to $\Field v$.

A \emph{symplectic basis} of $(\V,\omega)$ is a pair $(u,v)\in \V\times \V$ such that $\omega(u,v)=1$.
Denote by $S_{u,v}\in \GL(\V)$ the unique element of order $4$ sending $u$ to $v$, and by $\H_\omega \subset \SL(\V)\cap \Sl(\V)$ the set of such $S_{u,v}$.
Let $\Cone_\omega$ be the set of elements $p\in \Cone$ such that there exists $S_{u,v}\in \H_\omega$ for which the scalar product $\langle p , S_{u,v} \rangle$ is a sum of squares $x^2+y^2$ of elements $x,y\in \Field$ (in which case this is true for all $S_{u,v} \in \H_\omega$).

\begin{figure}[h]
    \centering
    \scalebox{.5}{\input{images/tikz/lotus-quadratic-cone-projection}}
    \caption*{The quadratic map $\psi \colon \Q^2 \to \Cone$, and the isomorphism $\P(\psi) \colon \P(\Q^2)\to \P(\Cone)$ of projective lines. Consult \cite[Chapter 2]{CLS_phdthesis_2022} for details about this construction.} 
\end{figure}

\begin{Lemma}
\label{Lem:param-cone-omega}
The quadratic map $\psi \colon \V \to \Cone$ has image $\Cone_\omega$ and is two-to-one outside the origin.
It intertwines the tautological action of $\SL(\V)$ on $\V$ with the restriction of its adjoint action on $\Cone_\omega$.
For all $u,v \in \V$ we have $2\langle \psi(u), \psi(v) \rangle = \omega(u,v)^2$.

We have $\P(\Cone_\omega)=\P(\Cone)$ and the map $\P(\psi)\colon \P(\V)=\P(\Cone_\omega)$ is an isomorphism of projective lines, inverse to the restriction of $\Psi$ to $\P(\Cone)$.
%
%
For distinct $u,v,x,y \in \Cone_\omega$ we have:
\begin{equation*}
    \label{cross-ratio-scalar}
    \tag{CRS}
    \bir(u,v,x,y)^2= 
    \tfrac{\langle u,v \rangle \cdot \langle x,y \rangle}{\langle u,y \rangle \cdot \langle x,v \rangle
    }.
\end{equation*}

\end{Lemma}

\begin{proof}
One may show that $\psi(\V)=\Cone_\omega$ and $\P(\Cone_\omega)=\P(\Cone)$ by computing in a symplectic basis of $\V$.
The last equality equals the cross-ratio of the four lines in $\V$ generated by preimages of $u,v,x,y$, computed in terms of the area form, as one can see in the following Figure.
\begin{figure}[h]
    \centering
    \scalebox{0.42}{
\begin{tikzpicture}[line cap=round,line join=round,>=triangle 45,x=2.0cm,y=2.0cm]
\clip(-4.1,-0.3) rectangle (4.1,1.3);

\draw[line width=1.5pt,color=forestgreen,-{Stealth[length=3.5mm,width=3mm]}] (0,0) -- (0,1);
\draw[line width=1.5pt,color=red,-{Stealth[length=3.5mm,width=3mm]}] (0,0) -- (1,1);
\draw[line width=1.5pt,color=blue,-{Stealth[length=3.5mm,width=3mm]}] (0,0) -- (-1,1);
\draw[line width=1.5pt,color=violet,-{Stealth[length=3.5mm,width=3mm]}] (0,0) -- (0.5,1);

\fill[color=blue,opacity=0.2](0,0) -- (-1,1) -- (0,1) -- cycle;
\fill[color=forestgreen,opacity=0.2](0,0) -- (0.5,1) -- (0,1) -- cycle;
\fill[color=red,opacity=0.2](0,0) -- (0.5,1) -- (1,1) -- cycle;

\begin{scriptsize}
\draw[color=blue,anchor=south] (-1,1) node {\Huge$\textcolor{blue}{\vec{x}}$};
\draw[color=forestgreen,anchor=south] (0,1) node {\Huge$\vec{v}$};
\draw[color=violet,anchor=south] (0.5,1) node {\Huge$\vec{y}$};
\draw[color=red,anchor=south] (1,1) node {\Huge$\vec{u}$};

\draw[fill=black] (0,0) circle (2pt);
\draw[color=black,anchor=north east] (0,0) node {\normalsize$O$};

\draw[color=black] (-3,0.5) node {\Huge$\displaystyle \frac{\textcolor{blue}{X}\textcolor{violet}{Y}}{\textcolor{blue}{X}\textcolor{forestgreen}{V}}=\frac{\omega(\textcolor{blue}{\vec{x}},\textcolor{violet}{\vec{y}})}{\omega(\textcolor{blue}{\vec{x}},\textcolor{forestgreen}{\vec{v}})}$};

\draw[color=black] (3,0.5) node {\Huge$\displaystyle \frac{\textcolor{red}{U}\textcolor{forestgreen}{V}}{\textcolor{red}{U}\textcolor{violet}{Y}}=\frac{\omega(\textcolor{red}{\vec{u}},\textcolor{forestgreen}{\vec{v}})}{\omega(\textcolor{red}{\vec{u}},\textcolor{violet}{\vec{y}})}$};
\end{scriptsize}
\end{tikzpicture}
}
    \label{fig:cross-ratio-det}
\end{figure}
\end{proof}

\subsection*{Ptolemy's theorem for ideal quadrilatareals}

We now apply the previous Lemma to show the following analogue of Ptolemy's theorem for quadrilaterals inscribed in the projective conic $\P(\Cone)$, which appeared in \cite[Proposition 2.6]{Penner_deco-teich-space-punct-surf_1987}.
It is better formulated if we fix a symplectic form $\omega$ on $\V$ and consider vectors in the subset $\Cone_\omega$ of the isotropic cone $\Cone$.
Since $\P(\Cone_\omega)=\P(\Cone)$ we may always lift a quadrilateral to such a quadruple, and any lift will satisfy the identity.

\begin{Proposition}
\label{Prop:Ptolemy-cone}
For distinct $u,v,x,y \in \Cone_\omega$, the following identity holds in an extention of $\Field$:

\begin{equation} \label{Ptolemy-ideal}\tag{IPS} 
    \sqrt{\langle u,v \rangle \cdot \langle x,y \rangle} =
    \sqrt{\langle u,y \rangle \cdot \langle x,v \rangle} +
    \sqrt{\langle u,x \rangle \cdot \langle v,y \rangle}
\end{equation}
\end{Proposition}

This formula is invariant under the action of $(\Field^\times)^2$ by individual dilatation of $u,v,x,y$, so we may suppose they lie on a conic section $\{p\in \Cone_\omega \mid \langle S_{u,v},p \rangle = 1\}$.

\begin{proof} 
This identity is equivalent, after dividing by the left hand side, to: \begin{equation*}
    \sqrt{\tfrac{\langle u,v \rangle \cdot \langle x,y \rangle}{\langle u,y \rangle \cdot \langle x,v \rangle}} + 
    \sqrt{\tfrac{\langle u,v \rangle \cdot \langle x,y \rangle}{\langle u,x \rangle \cdot \langle v,y \rangle}} = 1.
\end{equation*}
But this follows from the identity in the previous Lemma \ref{Lem:param-cone-omega} and the addition rule of cross-ratios $\bir(u,v,x,y)^{-1}+\bir(u,v,y,x)^{-1}=1$.
%
\end{proof}

\newpage

\section{The adjoint actions of \texorpdfstring{$\PGL_2(\Field)$}{PGL(2;K)} and \texorpdfstring{$\PSL_2(\Field)$}{PSL(2;K)} on \texorpdfstring{$\P(\Sl_2(\Field))$}{P(sl(2;K))}}
\label{Sec:adjoint-action}

%
%

\subsection*{The isomorphism \texorpdfstring{$\PGL(\V) \to \SO(\Sl(\V),\det)$}{PGL(2;K)=SO(2,1)}}

The left adjoint linear action of $\GL(\V)$ on $\Gl(\V)$ preserves the involution whence every structure which derives from it, such as the determinant form and the orthogonal decomposition $\Field \Id \oplus \Sl(\V)$. 
%
%
It also preserves the orientations of $\Sl(\V)$ defined for a basis as the class of its determinant in $\Field^\times/(\Field^\times)^2$.
Only the scalar matrices act trivially, and the maximal subspace on which the action is trivial equals $\Field\Id$. Therefore no information is lost after quotienting by these centers, and this yields a faithful representation $\PGL(\V) \to  \SO(\Sl(\V),\det)$ into the group of orientation preserving isometries of $(\Sl(\V),\det)$.


\begin{Proposition}
\label{Prop:iso-adjoint-action}
The adjoint action yields an isomorphism $\PGL(\V) \simeq \SO(\Sl(\V),\det)$.
\end{Proposition}

\begin{proof}
To prove surjectivity, we use a theorem of Cartan-Dieudonn\'e \cite{Dieudonne_geometrie-groupes_1971} stating that every isometry of a symmetric non-degenerate bilinear form over a $d$-dimensional $\Field$-vector space is a product of at most $d$ reflections.
In particular, an element of $\SO(\Sl_2(\Field),\det)$ is a product of at most $3$ reflections, but since it has determinant $1$ it is in fact a product of exactly two reflections. Thus we must express all products of two reflections as the conjugacy by some element.

If $q\in \Gl(\V)$ is not isotropic, that is $\det(q)\ne0$, then the orthogonal reflection $\sigma_q \in \End(\Gl(\V))$ of vector $q$ across $q^\perp$ is given by: 
\begin{equation*}
    \sigma_q(m) = m-2\tfrac{\langle q,m \rangle}{\langle q,q \rangle}\cdot q
\end{equation*}
Notice that the orthogonal reflection of vector $\Id$ across $\Sl(\V)$ equals $\sigma_{\Id} \colon m\mapsto -m^\#$.
The endomorphism $\mu_q \in \End(\Gl(\V))$ corresponding to left multiplication by $q$, left conjugates $\sigma_{\Id}$ to $\sigma_q$.
In formulae, we have $\mu_q \colon m \mapsto qm$ and $\sigma_{q}=\mu_q \circ \sigma_{\Id} \circ \mu_{q^{-1}}$.
Thus
\begin{equation*}
    \sigma_q(m)
    = -q(q^{-1}m)^\#
    = - \tfrac{qm^\#q}{\det(q)}.
\end{equation*}

Now restricting the attention to $\End(\Sl(\V))$, we notice that for $q,m\in \Sl(\V)$ this formula becomes $\sigma_q(m)=-qmq^{-1}$.
Hence for $p,q\in \Sl(\V)\setminus \Cone$ the reflection $\sigma_{p}\circ \sigma_{q} \in \SO(\Sl(\V),\det)$ coincides with the left adjoint action of $pq \in \GL(\V)$.
\end{proof}

The adjoint action commutes with the projectivization map $\Sl(\V) \to \P(\Sl(\V))$.
This realizes $\PGL(\V)$ as a subgroup inside the automorphism group $\PGL(\Sl(\V))$ of the projective plane $\P(\Sl(\V))$, namely the stabiliser of the non-degenerate conic $\P(\Cone)$.

The description of the actions of $\PGL(\V)$ and $\PSL(\V)$ on $\P(\Cone)$ follow from Lemma \ref{Lem:Param-Quadric-Two-Lines}.

\subsection*{Symmetric space of \texorpdfstring{$\PGL(\V)$}{PGL(V)}}

We call \emph{symmetry} of $\PGL(\V)$ an element of order two (since it maps to an orthogonal symmetry in $\SO(\Sl(\V),\det)$).
Those are represented by the elements in $\GL(\V) \cap \Sl(\V)=\Sl(\V)\setminus \Cone$, so the symmetries of $\PGL(\V)$ correspond by the projectivisation map to the complement $\P(\Sl(\V)\setminus \Cone)$ of the projective conic.
This is an open projective variety whose irreducible components over $\Field$ are indexed by the values of $\det \colon \GL(\V)\cap \Sl(\V)\to \Field^\times / (\Field^\times)^2$.
We call this variety $\P(\Sl(\V)\setminus \Cone)$ the \emph{symmetric space} of $\PGL(\V)$, in the spirit of \cite{Cartan_Lecons-3-reprint_1992}.

Hence the group $\PGL(\V)$ acts on its symmetric space $\P(\Sl(\V)\setminus \Cone)$ by the projectivised adjoint representation, and the elements of order two are the symmetries.
Since $s\in \GL(\V)\cap \Sl(\V)$ maps to an element of order two in $\SO(\Sl(\V),\det)$ which fixes the line $\Field s$, it acts like minus the identity on the orthogonal plane $s^\perp$, thus corresponds to the orthogonal symmetry across the line $\Field s$:
%
\begin{equation*}
    \forall x\in \Sl(\V) \:\colon\: \quad
    sxs^{-1} + x = 2\tfrac{\langle s,x\rangle }{\langle s,s\rangle}\cdot s
\end{equation*}
and we recognise from the proof of Proposition \ref{Prop:iso-adjoint-action}, the expression for the composition of reflections $\sigma_s \circ \sigma_\Id \in \SO(\Gl(\V),\det)$ restricted to $\Sl(\V)$.

\subsection*{Action of \texorpdfstring{$\PGL(\V)$}{PGL(V)} on \texorpdfstring{$\P(\Sl(\V)\setminus \Cone)$}{P(sl(V-X)}}

Let us begin with another corollary to Lemma \ref{Lem:param-cone-omega}, which describes the adjoint action of an element $C\in \PSL(\V)$ on $\Sl(\V)$.

\begin{Corollary}
\label{Cor:action-A-perp}
%
For $C\in \SL(\V)$, the adjoint action of $C$ on $\Sl(\V)$ restricted to the plane $(\pr C)^\perp$ is equivalent over $\Field$ to the tautological action of $C^2$ on $\V$.
\end{Corollary}

\begin{Proposition}
\label{Prop:Conj-a-b-K[a,b]}
Consider distinct $\Aa,\Bb\in \Sl(\V)$ with determinant $d \ne 0$ and $\bir(\Aa,\Bb) \ne \infty$. 

The quadratic subalgebra $\Field[\{\Aa,\Bb\}]$ of $\Gl(\V)$ contains a unique $M\in \GL(\V)$ with $\Tr(M) = 2$ which conjugates $\Aa$ to $\Bb$. It is given by:
\begin{equation*}
    M 
    = \Id+ \tfrac{\bir(\Aa,\Bb)}{2d}\cdot \{\Aa,\Bb\}
    = \frac{(d+\langle \Aa,\Bb\rangle)\Id+ \{\Aa,\Bb\}}{d+\langle \Aa,\Bb\rangle}
    \qquad
    \mathrm{and}
    \quad \det(M) = \bir(\Aa,\Bb).
\end{equation*}
\end{Proposition}

\begin{proof}
First suppose $d=1$, so $\Aa,\Bb\in \H$.
For $x\in \Field$, set $M=\Id+x\{\Aa,\Bb\}$. We have $M\Aa=\Bb M \iff \Aa+x\{\Aa,\Bb\}\Aa=\Bb+x\Bb\{\Aa,\Bb\}$.
But $\{\Aa,\Bb\}\Aa=\frac{1}{2}(\Aa\Bb\Aa+\Bb)$, and since $\Aa\in \H$ acts like a symmetry across $\Field \Aa$, we have $\Aa\Bb\Aa=-\Aa\Bb\Aa^{-1}=\Bb-2\langle \Aa,\Bb\rangle \Aa$ so $\{\Aa,\Bb\}\Aa=\Bb-\langle \Aa,\Bb\rangle \Aa$.
Similarly $\Bb\{\Aa,\Bb\}=\Aa-\langle \Aa,\Bb\rangle \Bb$.
Thus $M\Aa=\Bb M \iff (\Aa-\Bb)(1-x(1+\langle \Aa,\Bb \rangle))=0 \iff 1=x(1+\langle \Aa,\Bb \rangle)$ since $\Aa-\Bb\ne 0$.

Now suppose $\Aa,\Bb\in \Sl(\V)$ have the same determinant $d\ne 0$. 
Divide them by $\sqrt{d}$, which may live in a quadratic extension $\Field'$ of $\Field$, to get $\Aa',\Bb'\in \H$ as before with $\bir(\Aa, \Bb) = \bir(\Aa', \Bb')$. Since $\{\Aa,\Bb\}/d=\{\Aa',\Bb'\}$ we have $\Field'[\{\Aa,\Bb\}]=\Field'[\{\Aa',\Bb'\}]$, and for $M\in \Field'[\{\Aa,\Bb\}]^\times$ an invertible element of this quadratic algebra, we have $M\Aa=\Bb M \iff M\Aa'=\Bb'M$ which completes the proof. 
\end{proof}

\begin{figure}[h]
    \centering
    \scalebox{0.5}{
\begin{tikzpicture}[x=4cm,y=4cm,z=2cm]

\path (0.,0.,0.) coordinate (O)
(0,-1,0) coordinate (O')
(1,0.,0.) coordinate (A)
(1,-2.,0.) coordinate (A')
(-1,0.,0.) coordinate (B)
(-1,-2,0.) coordinate (B')
(1./1.3,0,0.5/1.3) coordinate (D)
(-0.5,0,-0.5) coordinate (E)
(1./1.3-0.43,-2,0.5/1.3) coordinate (D')
(-0.5-0.43,-2,-0.5) coordinate (E');



\draw[name path=ElliH,line width=2.pt] (0,0,0) ellipse (1 and 0.15);
\draw[name path=ElliBBack,line width=2.pt,dashed] (1,-2,0) arc(0:180:1 and 0.15);
\draw[line width=2.pt,name path=ElliBFront] (1,-2,0) arc(0:-180:1 and 0.15);

\draw[name path=cone1,line width=1.5pt, color=black] (B) -- (A') ;
\draw[name path=cone2,line width=1.5pt, color=black] (A) -- (B') ;

\draw[name path=DE,line width=1.5pt, color=grey] (D) -- (E) ;
\draw[name path=D'E',line width=1.5pt, color=grey,dashed] (D') -- (E') ;
\draw[name path=EE',line width=1.5pt, color=grey] (E) -- (E') ;
\draw[name path=DD',line width=1.5pt, color=grey,dashed] (D') -- (D) ;
\path [name intersections={of=DE and ElliH,by={h1,h2}}];
\path [name intersections={of=D'E' and ElliBBack,by={bb}}];
\path [name intersections={of=D'E' and ElliBFront,by={bf}}];

\draw[line width=1.5pt, color=grey,dashed] (h1) -- (bf) ;
\draw[line width=1.5pt, color=grey,dashed] (h2) -- (bb) ;
\draw[line width=1.5pt, color=grey] (O') -- (bf) ;
\draw[line width=1.5pt, color=grey] (h2) -- (O') ;

\path [name intersections={of=D'E' and ElliBFront,by={i1}}];
\path [name intersections={of=D'E' and ElliBBack,by={i2}}];
\draw[line width=1.5pt, color=grey] (E') -- (i1) ;

\path [name intersections={of=DD' and cone1,by={d1}}];
\path [name intersections={of=DD' and ElliH,by={d2,d3}}];
\draw[line width=1.5pt, color=grey] (d1) -- (d3) ;
\draw[line width=1.5pt, color=grey] (d2) -- (D) ;


\coordinate (a) at ($(h1)!0.3!(h2)$);
\coordinate (b) at ($(h1)!0.7!(h2)$);
\coordinate (c) at ($(a)!0.6!(b)$);

\draw[line width=2pt, color=violet] (a) -- (b);
\draw[line width=2pt, color=violet, -{Stealth[length=4mm,width=3mm]}] (a) -- (c);

\node[draw,color=red,circle,inner sep=2pt, fill=red, label=-90:\Huge$\textcolor{red}{a}$] at (a) {};
\node[draw,color=blue,circle,inner sep=2pt, fill=blue, label=90:\Huge$\textcolor{blue}{b}$] at (b) {};

\draw[line width=2pt, color=violet,-{Stealth[length=4mm,width=3mm]}] (O') -- (0.7,-0.5,0);
\draw[line width=2pt, color=violet,anchor=south] (0.55,-.95,0) node {\Huge$\{\textcolor{red}{a},\textcolor{blue}{b}\}$};

\coordinate (h3) at ($(h2)!0.5!(b)$);
\coordinate (h4) at ($(h3)!0.3!(O')$);

\coordinate (h5) at ($(h1)!0.5!(a)$);
\coordinate (h6) at ($(h5)!0.3!(O')$);

\draw [color=violet,line width=2.pt,line cap =round,-{Stealth[length=4mm,width=3mm]}] (h6) to[out=-120, in=-70,looseness=3]  (h4);

\coordinate (i3) at ($(i2)!0.2!(i1)$);
\coordinate (i4) at ($(i3)!0.3!(O')$);

\coordinate (i5) at ($(i1)!0.2!(i2)$);
\coordinate (i6) at ($(i5)!0.3!(O')$);

\draw [color=violet,line width=2.pt,line cap =round,{Stealth[length=4mm,width=3mm]}-] (i4) to[out=100, in=60,looseness=3]  (i6);







\end{tikzpicture}
}
    \hfill
    \scalebox{0.5}{

\begin{tikzpicture}[line cap=round,line join=round,>=triangle 45,x=2.0cm,y=2.0cm]
\clip(-2.1,-2.1) rectangle (3.6,2.1);
\draw [line width=1.2pt] (0,0) circle (4cm);
\draw [dash pattern=on 2pt off 2pt,color=red,domain=-2.1:3.6] plot(\x,{(-4.61--1.78*\x)/-1.46});
\draw [dash pattern=on 2pt off 2pt,color=blue,domain=-2.1:3.6] plot(\x,{(--4.61-1.19*\x)/-1.97});
\draw [line width=1.2pt,color=violet,domain=-2.1:3.6] plot(\x,{(--3.96-2.98*\x)/-0.51});
\draw [line width=1.6pt,color=violet] (1.38,0.29)-- (1.18,-0.89);
\draw [line width=1.6pt,color=violet] (1.27,-0.34) -- (1.23,-0.29);
\draw [line width=1.6pt,color=violet] (1.27,-0.34) -- (1.33,-0.31);
\draw [line width=1.6pt,color=violet] (1.29,-0.26) -- (1.24,-0.21);
\draw [line width=1.6pt,color=violet] (1.29,-0.26) -- (1.35,-0.22);
\draw [line width=1.6pt,color=violet] (1.26,-0.43) -- (1.21,-0.38);
\draw [line width=1.6pt,color=violet] (1.26,-0.43) -- (1.32,-0.39);
\draw [line width=1.6pt,color=violet] (1.55,1.27)-- (1.38,0.29);
\draw [line width=1.6pt,color=violet] (1.46,0.74) -- (1.41,0.79);
\draw [line width=1.6pt,color=violet] (1.46,0.74) -- (1.52,0.77);
\draw [line width=1.6pt,color=violet] (1.47,0.82) -- (1.43,0.87);
\draw [line width=1.6pt,color=violet] (1.47,0.82) -- (1.53,0.85);
\draw [line width=1.6pt,color=violet] (1.44,0.65) -- (1.4,0.7);
\draw [line width=1.6pt,color=violet] (1.44,0.65) -- (1.5,0.68);
\draw [line width=1.6pt,color=violet] (1.18,-0.89)-- (1.04,-1.71);
\draw [line width=1.6pt,color=violet] (1.1,-1.34) -- (1.05,-1.29);
\draw [line width=1.6pt,color=violet] (1.1,-1.34) -- (1.16,-1.31);
\draw [line width=1.6pt,color=violet] (1.11,-1.26) -- (1.07,-1.21);
\draw [line width=1.6pt,color=violet] (1.11,-1.26) -- (1.17,-1.22);
\draw [line width=1.6pt,color=violet] (1.08,-1.43) -- (1.04,-1.38);
\draw [line width=1.6pt,color=violet] (1.08,-1.43) -- (1.14,-1.39);
\begin{scriptsize}
\fill [color=violet] (3.01,-0.52) circle (2.5pt);
\draw[color=violet] (3.1,-0.1) node {\Huge$\{\textcolor{red}{a},\textcolor{blue}{b}\}$};
\fill [color=red] (1.55,1.27) circle (1.5pt);
\fill [color=blue] (1.04,-1.71) circle (1.5pt);
\fill [color=red] (1.38,0.29) circle (2.5pt);
\draw[color=red] (0.98,0.36) node {\Huge$a$};
\fill [color=blue] (1.18,-0.89) circle (2.5pt);
\draw[color=blue] (0.81,-0.79) node {\Huge$b$};
\end{scriptsize}
\end{tikzpicture}
}
    \caption*{The one parameter group generated by $\{a,b\}$, which is contained in $\Span(\Id,\{a,b\})$, acts by translation along the line $(a,b)$.}
    \label{fig:sqrt(ab)}
\end{figure}

\begin{Corollary}
The group $\PGL(\V)$ acts transitively on each level set of the determinant in $\Sl(\V)\setminus \Cone$, and therefore on each irreducible component of its symmetric space $\P(\Sl(\V) \setminus \Cone)$.

The semi-simple conjugacy classes in $\PGL(\V)$ are classified by the value of $(\Tr A)^2/\disc(A)\in \Field$, which is $0$ for the class of involutions.
\end{Corollary}

%

Recall that if $\det(\Aa)\ne 0$ then the stabiliser of $\Aa\in \Sl_2(\Field)\setminus \Cone$ in $\Gl_2(\Field)$ is reduced to the quadratic subalgebra $\Field[\Aa]$. This implies the following Corollary, which will have arithmetic applications bearing to the genus of quadratic forms.

\begin{Corollary}
\label{Cor:bir-mod-norms}
Consider distinct $\Aa,\Bb\in \Sl(\V)$ with determinant $-\delta \ne 0$ and $\bir(\Aa,\Bb) \ne \infty$.

The matrices $M\in \PGL_2(\V)$ conjugating $\Aa$ to $\Bb$ have a well defined determinant in the quotient $\Field^\times/\Norm_\Field(\Field[\sqrt{\delta}]^\times)$, and its is equal to the class of $\bir(\Aa,\Bb)$.
\end{Corollary}



Let us now describe the structure of the orbits for the adjoint action of $\PSL(\V)$ on the non-zero level sets $\{\det = -\delta\}\subset \Sl(\V)$, which is the main Theorem in \cite[Chapter 1]{CLS_phdthesis_2022}.

\begin{Theorem}
\label{Thm:Conj-a-b-SL(V)}
Let $\Aa,\Bb\in \Sl(\V)$ have determinant $-\delta \ne 0$ and cross-ratio $\bir(\Aa,\Bb)=4\chi \notin \{1,\infty\}$.
The elements $C\in \SL(\V)$ such that $C\Aa C^{-1}= \Bb$ are parametrized by the Pell-Fermat conic:
\begin{equation*}
    (x,y)\in \Field\times \Field 
    \: \colon\: \quad
    x^2-\delta y^2 = \chi
\end{equation*}
\begin{equation*}
    C(x,y) = x(\Id+\Bb\Aa^{-1})+y(\Aa+\Bb)
\end{equation*}

In particular, $\Aa$ and $\Bb$ are conjugate by an element $C(x,y)\in \SL(\V)$ if and only if $\bir(\Aa,\Bb)$ belongs to the subgroup of norms $\Norm_\Field \Field[\sqrt{\delta}]\subset  \Field^\times$ of the quadratic extension, and by an element $C(x,0)\in \SL(\V) \cap \Field[\{\Aa,\Bb\}]$ if and only if $\bir(\Aa,\Bb)$ belongs to the subgroup of squares $(\Field^\times)^2 \subset  \Field^\times$.
\end{Theorem}

\begin{proof}
Suppose first that $\Aa,\Bb\in \H$.
Since $\det\{\Aa,\Bb\}\ne 0$ the elements $\Id,\Aa,\Bb,\{\Aa,\Bb\}$ form a basis of $\Gl(\V)$.
%
Let $C\in \Gl(\V)$ be decomposed as $C=t\Id+x\{\Aa,\Bb\}+y\Aa+z\Bb$ for $t,x,y,z\in \Field$.
The condition $C\Aa=\Bb C$ can be rewritten using $\Aa^2=-\Id=\Bb^2$ as well as $\{\Aa,\Bb\}\Aa=\Bb-\langle \Aa,\Bb\rangle \Aa$ and $\Bb\{\Aa,\Bb\}=\Aa-\langle \Aa,\Bb\rangle \Bb$.
After grouping terms we find $C\Aa=\Bb C \iff \left(t-x(1+\langle \Aa,\Bb\rangle)\right)\cdot (\Aa-\Bb)+(z-y)\cdot (\Id+\Bb\Aa)=0$.
%
But $\Aa-\Bb \in \Span(\Aa,\Bb)\setminus\{0\}$ and $\Id+\Bb\Aa=(1-\langle \Aa,\Bb \rangle)\Id - \{\Aa,\Bb\}\in \Field[\{\Aa,\Bb\}]\setminus\{0\}$ so by orthogonality of the planes $\Field[\{\Aa,\Bb\}]$ and $\Span(\Aa,\Bb)$ we have $C\Aa=\Bb C \iff t = x(1+\langle \Aa,\Bb\rangle) \:\&\: y=z$. 

Now for $x,y\in \Field$ the determinant of $C=x(\Id-\Bb\Aa)+y(\Aa+\Bb)$ can be computed using the orthogonality of $\Span(\Id,\Aa\Bb)$ and $\Span(\Aa,\Bb)$ and the hypothesis $\Aa,\Bb\in \H$:
\begin{equation*}
    \det(C)=(x^2+y^2)\cdot(2+2\langle \Aa,\Bb\rangle) 
    = \tfrac{4(x^2+y^2)}{\bir(\Aa,\Bb)}
\end{equation*}
so $C\in \SL(\V) \iff \bir(\Aa,\Bb)=(2x)^2+(2y)^2$.
This proves the Lemma for $\delta = -1$.

Finally, let us reduce the general case $\delta \ne 0$ to the previous one.
The points $\Aa'=\Aa/\sqrt{d}$ and $\Bb'=\Bb/\sqrt{d}$ satisfy $\bir(\Aa,\Bb)=\bir(\Aa',\Bb')$ and the endomorphisms $C\in \SL(\V\otimes \Field')$ with $\Field'=\Field[\sqrt{d}]$ conjugating $\Aa$ to $\Bb$ are the same as those conjugating $\Aa'$ to $\Bb'$.
We just showed that those elements $C$ correspond to the pairs $(x,y)\in \Field'\times \Field'$ such that $\bir(\Aa,\Bb)=4(x'^2+y'^2)$ by the formula $C = x'(\Id + \Bb'\Aa'^{-1}) + y' (\Aa'+\Bb')$. 
Setting $x=x'$ and $y = y'\tfrac{1}{\sqrt{d}}$ which satisfy $4(x^2+dy^2) = \bir(\Aa,\Bb)$, we may rewrite $C = x(\Id + \Bb\Aa^{-1}) + y (\Aa+\Bb)$. 
But recall that $a$ and $b$ have coefficients $\Field$ and that $\Aa+\Bb\ne 0$ is orthogonal to $\Id+\Bb\Aa^{-1}\ne 0$. Hence $C$ is has coefficients in $\Field$ if and only if $x,y\in \Field$.
\end{proof}

\begin{Remark}
Theorem \ref{Thm:Conj-a-b-SL(V)} holds with coefficients restricted to any subring of $\Field$ containing $1/2$.
For instance $\Z[1/2] \subset \Q$, or $\Z_p \subset \Q_p$ for odd prime $p$.
\end{Remark}

\begin{Remark}
The problem of conjugating $\Aa,\Bb \in \Sl(\V)\setminus \Cone$ by $C\in \SL(\V)$ can be formulated as the search for fixed points under the transformation $C \mapsto \Aa C \Bb^{-1}$. Let us be more precise.

The group $\Gl(\V)^\times \times \Gl(\V)^\times$ acts linearly on $\Gl(\V)$ by $(A,B)\cdot C = A C B^{-1}$, preserving the isotropic cone $\{\det=0\}$ and its complement $\Gl(\V)^\times$.
The $C\in \Gl(\V)^\times$ conjugating $A,B\in \Gl(\V)^\times$ correspond to the fixed points of $(A,B)$ under this action.
If they exist, then $(A,B)$ must belong to the subgroup of pairs with $\det(A)=\det(B)$, preserving the level sets of $\det$, namely the stabiliser of $\SL(\V)$.

The proof of Theorem \ref{Thm:Conj-a-b-SL(V)} can be recast as the description of the linear action $(\Aa,\Bb)$ on $\Gl(\V)$, which under the assumptions $\det(\Aa)=-\delta=\det(\Bb)$ and $\langle a, b \rangle = \kappa  \ne \pm \delta $ ensuring that $(\Id, \{\Aa,\Bb\},\Aa,\Bb)$ forms a basis, is given by the following matrix where $c=-\kappa/\delta$ is the cosine:
\begin{equation*}
    L_{\Aa} = 
    \begin{psmallmatrix}
    0 & 0 & \delta & -\kappa \\
    0 & 0 & 0 & 1 \\
    1 & \kappa & 0 & 0 \\
    0 & \delta & 0 & 0
    \end{psmallmatrix}
    \quad
    R_{\Bb^{-1}} = \tfrac{1}{\delta}
    \begin{psmallmatrix}
    0 & 0 & -\kappa & \delta \\
    0 & 0 & 1 & 0 \\
    0 & \delta & 0 & 0 \\
    1 & \kappa & 0 & 0
    \end{psmallmatrix}
    \qquad
    L_{\Aa}R_{\Bb^{-1}} = \tfrac{1}{\delta}
    \begin{psmallmatrix}
    -\kappa & \delta^2-\kappa^2 & 0 & 0 & \\
    1 & \kappa & 0 & 0\\
    0 & 0 & 0 & \delta \\
    0 & 0 & \delta & 0
    \end{psmallmatrix}
    = 
    \begin{psmallmatrix}
    c & (1-c^2)\delta & 0 & 0 \\
    1/\delta & -c & 0 & 0\\
    0 & 0 & 0 & 1 \\
    0 & 0 & 1 & 0
    \end{psmallmatrix}
\end{equation*}

Its characteristic polynomial is $(X^2+1)^2 = (X-1)^2(X+1)^2$ and some eigenvectors for the eigenvalues $+1$ and $-1$ are $C(x,y)^{\#}=(x(1+c),x/\delta,-y,-y)$ and $(x(1+c),-x/\delta,y,-y)$.
\end{Remark}

\begin{Remark}
We have 
$\Tr C(x,y) = 2x(1-\langle \Bb,\Aa \rangle/\delta) = 2x(1+\cos(\Aa,\Bb)) = 4x/\bir(\Aa,\Bb) = x/\chi$.
\end{Remark}

\begin{Remark}[Square roots of $\Bb\Aa^{-1}$]
For $y=0$ we recover the unique multiples of $M$ satisfying the conditions in Proposition \ref{Prop:Conj-a-b-K[a,b]} which belong to $\SL(\V\otimes \Field')$ where $\Field'=\Field[\sqrt{\bir(\Aa,\Bb)}]$, namely:
\begin{equation*}
    \pm C =\tfrac{1}{2}\sqrt{\bir(\Aa,\Bb)} \left(\Id+\Bb\Aa^{-1}\right).
\end{equation*}

These are the unique square roots of the product of symmetries $\Bb\Aa^{-1}$ in the extended quadratic algebra $\Field'[\{\Aa,\Bb\}]$ since one may compute that $-\delta C^2=\langle \Aa,\Bb \rangle + \{\Aa,\Bb\} = -\Bb\Aa$, thus
\begin{equation*}
    C^2 =\Bb\Aa^{-1}.
\end{equation*}

%
%
%
\end{Remark}

Our last proposition completes the description for the $\PSL(\V)$-orbits of pairs $\Aa,\Bb\in \Sl(\V)\setminus \Cone$. Two such pairs are conjugate if and only if the obvious conditions on the scalar products hold, together with the Pell-Fermat conditions on the cross ratios given by Theorem \ref{Thm:Conj-a-b-SL(V)}.

\begin{Proposition}
\label{Prop:Stab-a-transitive}
The stabiliser $\SL(\V)\cap \Field[\Aa]$ of $\Aa\in \Sl(\V)\setminus \Cone$ acts transitively on the set of elements in $\Sl(\V)$ with a given determinant and scalar product with $\Aa$.
\end{Proposition}
\begin{proof}
For $C=t\Id+u\Aa\in \SL(\V)$ and $\Bb,\Bb'\in \Sl(\V)$, we have $Cb=b'C$ if and only if $\langle \Aa,\Bb\rangle = \langle \Aa,\Bb'\rangle$ \& $t(\Bb-\Bb')+u\{\Aa,\Bb+\Bb'\}=0$. This last condition amounts to the colinearity of $(\Bb-\Bb')$ with $\{\Aa,\Bb+\Bb'\}$, that is $(\Bb-\Bb')\perp \Aa$ \& $(\Bb-\Bb')\perp (\Bb+\Bb')$, so the claim follows.
(If $\Bb\ne \Aa$ then $\pm C$ is unique.)
\end{proof}

Consider the action of $\PSL(\V)$ by conjugacy on itself.
Let us say that $A,B\in \PSL(\V)$ are of the same type if $\disc(A)\equiv \disc(B)\bmod{(\Field^\times)^2}$, in which case we may 
define $\bir(A,B)$.

\begin{Corollary}
%
A pairs of semi-simple elements $A_1,A_2$ of the same type is conjugate to another pair of semi-simple elements $B_1,B_2$ of the same type if and only if we have $\bir(A_1,A_2)=\bir(B_1,B_2)$ as well as $\disc(A_i)=\Delta_i=\disc(B_i)$ and $\bir(A_i,B_i)\equiv 1 \bmod{\Norm_\Field \Field(\sqrt{\Delta_i})^\times}$.
\end{Corollary}

\subsection*{Examples over finite Fields}

Let us fix a basis of $\V$ to identify it with $\Field^2$, and first deduce a canonical basis for $\Gl(\V)=\Gl_2(\Field)$. Its elements together with their opposites forms the dihedral group of order $8$ which acts faithfully on the square whose vertices have coordinates $\pm 1$:
\begin{equation*}
\Id = 
\begin{psmallmatrix}
1 & 0 \\ 
0 & 1
\end{psmallmatrix}
\qquad
S = 
\begin{psmallmatrix}
0 & -1 \\
1 & 0
\end{psmallmatrix}
\qquad
J = 
\begin{psmallmatrix}
0 & 1 \\
1 & 0
\end{psmallmatrix}
\qquad
K = 
\begin{psmallmatrix}
-1 & 0 \\
0 & 1
\end{psmallmatrix}
\end{equation*}
This orthogonal basis presents $\Gl_2(\Field)$ as a quaternion algebra, with product given by $J^2=K^2=\Id$ and $KJ=-JK=S$.
Its elements are defined over all $\Field$ since their coordinates belong to $\{-1,0,1\}$, in particular their properties are invariant under extension of scalars.

\begin{Example}[Conjugating elements with their opposites]
The pairs $(K,J)$ and $(-K,-J)$ are conjugate by $S\in \SL_2(\Field)$. The elements $S$ and $-S$ are conjugate by $K,J\in \GL_2(\Field)$ but according to Proposition \ref{Prop:Stab_GL_sl-X}, they are conjugate in $\SL_2(\Field)$ if and only if $-1$ is a sum of squares in $\Field$.
\end{Example}

\begin{Example}[Conjugating $J$ and $K$]
Let us show that $J$ and $K$ are conjugate by $\SL_2(\Field)$ but are conjugate by $\Field[\{J,K\}]$ only when $2\in (\Field^\times)^2$.

The $M\in \Gl_2(\Field)$ satisfying $MJ=KM$ are $M(p,q)=
\begin{psmallmatrix} 
p & -p \\ q & q \end{psmallmatrix}$ for $p,q\in \Field$, and $\det(M)=2pq$. 
The only ones in $\Field[\{J,K\}]=\Field[S]$ are those for which $p=q$.

Applying Theorem \ref{Prop:Conj-a-b-K[a,b]} to $J,K$ of determinant $-1$ with $\langle J,K\rangle = 0$ and $\{J,K\}=-S$ yields 
\begin{equation*}
    C(x,y)  = \tfrac{1}{\sqrt{2pq}}M(p,q)
    = \tfrac{1}{\sqrt{2pq}} \left[\tfrac{p+q}{2}(\Id + S) + \tfrac{q-p}{2}(J + K) \right] \in \SL_2(\Field[\sqrt{2pq}]).
\end{equation*}
%
Therefore $J$ and $K$ are always conjugate by $\SL_2(\Field)$ by choosing for instance $p=2q$, but they are conjugate in $\Field[\{J,K\}]$ only when $2\in (\Field^\times)^2$.
\end{Example}


\begin{Example}
Let us describe the orbits for the adjoint action of $\PSL_2(\F_3)$ on the non-zero level sets of $(\Sl_2(\F_3),\det)$, given by $1\colon \{\pm S, \pm (J \pm K)\}$ and $-1\colon \{\pm J, \pm K, \pm S \pm J \pm K\}$.
%
%

For $\Aa,\Bb\in \Sl_2(\F_3)$ with $\det=-\delta \ne 0$, the hypothesis $\bir(\Aa,\Bb)\notin\{1,\infty\}\iff\cos(\Aa,\Bb)^2 \ne 1$ becomes $\langle \Aa,\Bb \rangle = 0$, and leads to the equation $x^2-\delta y^2 = -1$ which always has a solution.
Hence Theorem \ref{Thm:Conj-a-b-SL(V)} says that any two orthogonal elements of the same non-zero determinant are conjugate.

We find for instance that $S$ is conjugated to $J+K$ by the element $1+S-K$.
Moreover as $-1$ is a sum of squares in $\F_3$ we know that $S$ and $-S$ are conjugate, for instance by $J+K\in \SL(\V)$.
From this we deduce that the level set $1$ consists of a single orbit with $6$ elements.

The quaternionic group $\pm \{\Id,S,J,K\}$ assembles the elements $\pm S \pm J \pm K$ in two orbits given by the product of their coefficients $\bmod{3}$.
Moreover we saw that $J,-J,K,-K$ are all conjugate.
An exhaustive search confirms that the $12$ element level set $\{\det = -1\}$ is partitioned into those $3$ orbits with $4$ elements, of which two form the vertices of tetrahedra whose edges intersect along the third forming the vertices of a square.
Using Proposition \ref{Prop:Stab-a-transitive}, we find that the action of  $\PSL_2(\F_3)$ on this picture recovers the isomorphism with the tetrahedral group $\triangle(2,3,3)\simeq \mathfrak{A}_4$.
%

One may also partition $\Cone\setminus\{0\}$ in two orbits forming the vertices of tetrahedra, whose edges correspond to level sets of the scalar product.

\begin{figure}[h]
    \centering
    \includegraphics[width=0.32\textwidth]{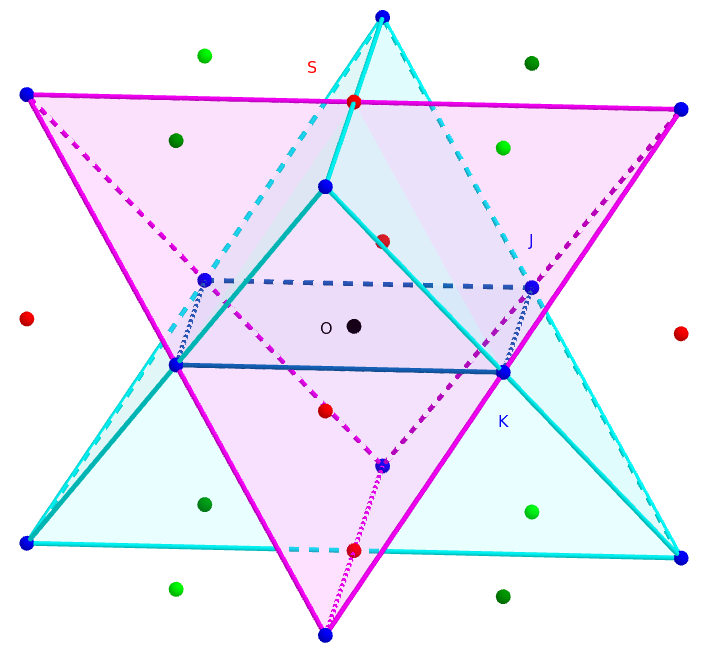}
    \hfill
    \includegraphics[width=0.31\textwidth]{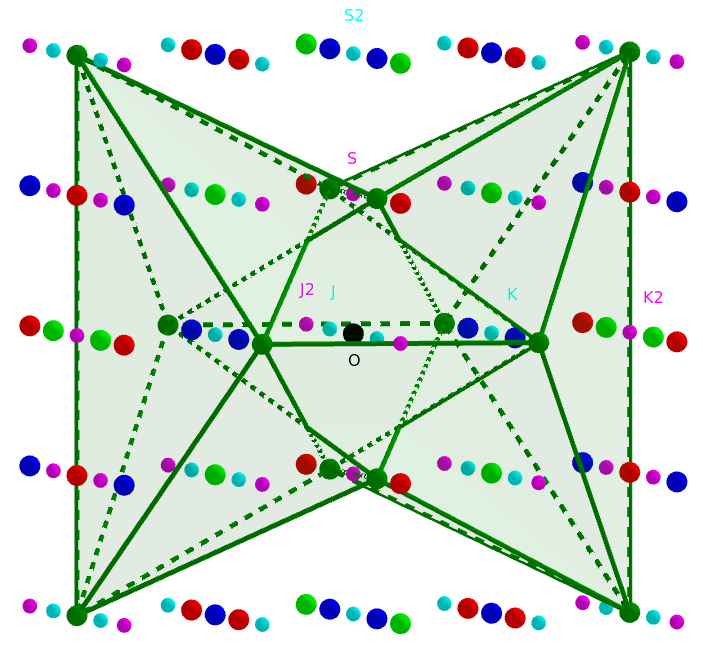}
    \hfill
    \includegraphics[width=0.31\textwidth]{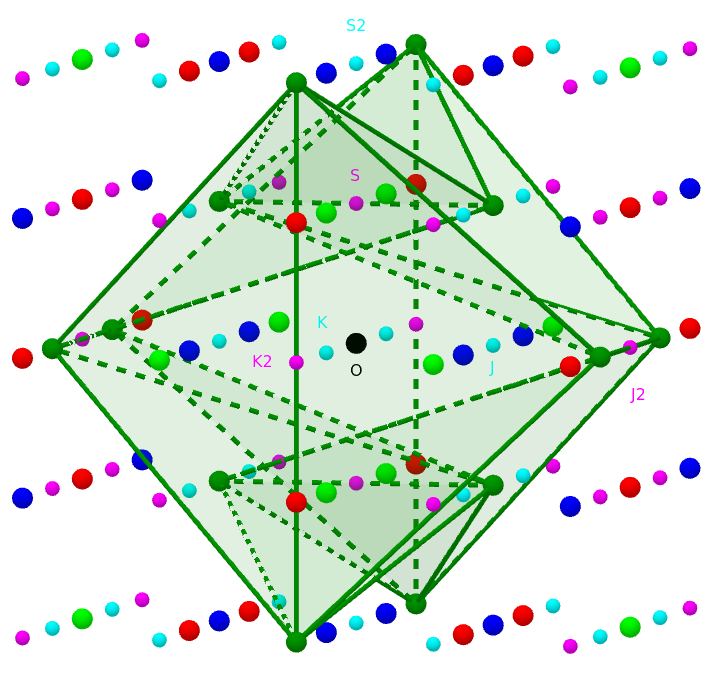}
    \caption*{$\Sl_2(\F_3)$: the \textcolor{green!50!black}{isotropic cone $\Cone$}, the orbit \textcolor{red}{$\{\det = +1\}$} and the partition of \textcolor{blue}{$\{\det = -1\}$} into $3$ orbits: the vertices of each tetrahedron and the intersections of their edges. $\Sl_2(\F_5)$ : the orbits $\det=$ \textcolor{red}{$-2$}, 
    \textcolor{magenta}{$-1$},
    \textcolor{cyan}{$1$},
    \textcolor{blue}{$2$} and an \textcolor{black!25!green}{icosahedral orbit} in $\Cone\setminus \{0\}$.}
\end{figure}

\end{Example}

\begin{Example}
In $\Sl_2(\F_5)$ the level sets of $\det$ for the values $1,2$ are given in the following table. Multiplication by $2$ yields a bijection between level sets with opposite signs, and the $25$ remaining elements belong to $\Cone$. 
Both forms $x^2+y^2$ and $x^2+2y^2$ represent all elements in $\F_5$, and one may use Theorem \ref{Thm:Conj-a-b-SL(V)} to show that the action of $\PSL_2(\F_5)$ is transitive on each non-zero level set of $\det$.

\begin{table}[h]
\centering
\begin{tabular}{|c|c|c|}
\hline
   $\det$
 & $\{kK+sS+jJ\mid k,s,j \in [-2,2]\}\subset \Sl_2(\F_5)$
 & Cardinal
 \\
 \hline
   $0$
 & \textcolor{green!50!black}{$2(\pm K \pm S),\; \pm K \pm 2J,\; \pm S \pm J$},\; $0,\;$ \textcolor{green}{$2(\pm K \pm S,\; \pm 2K \pm J,\; 2(\pm S \pm J)$}
 & $12+1+12$
 \\
\hline
   $1$
 & \textcolor{magenta}{$\pm S,\; \pm 2K,\; \pm 2J,\; \pm K\pm S\pm 2J,\; \pm 2K \pm S\pm J,\; \pm 2K \pm 2S\pm 2J$}
 & $30$
 \\
\hline
   $2$
 & \textcolor{red}{$\pm 2K\pm S,\; \pm S \pm 2J,\; \pm 2J\pm 2K,\; \pm J\pm 2S\pm K$}
 & $20$
 \\
\hline
\end{tabular}
\end{table}

The orbits $\{\det=\pm 2\}$ form the vertices of two octahedra whose edges correspond to the orbits $\{\det = \pm 1\}$, and Proposition \ref{Prop:Stab-a-transitive} shows that the action of $\PSL_2(\F_5)$ on this picture recovers the isomorphism with the icosahedral group $\triangle(2,3,5)\simeq \mathfrak{A}_5$.
The set $\Cone\setminus\{0\}$ partitions in two orbits forming the vertices of icosahedra, whose edges can be defined using level sets of the scalar product. 
\end{Example}

\newpage

\section{Applications to binary quadratic forms}
\label{Sec:Qfb}

Let $\Qfb(\V)$ be the space of quadratic forms $Q\colon \V \to \Field$. After choosing a basis of $\V$, those amount to homogeneous polynomials in two ordered variables with coefficients in $\Field$.

\subsection*{Isomorphism \texorpdfstring{$\Qfb(\V) \simeq \Sl(\V)$}{Qfb(V)=sl(V)}}
One may polarise $Q \in \Qfb(\V)$ with respect to any non degenerate bilinear form on the plane $\V$, and one usually learns this for some euclidean scalar product, but we may also use a symplectic form: there exists a unique $\Qq\in \Sl(\V)$ such that $Q(v)=\omega(v,\Qq v)$.
If we fix a basis $\V=\Field^2$ and $\omega = \det$ we have the formula:
\begin{equation*}
    Q=lx^2+mxy+ry^2 
    \in \Qfb(\Field^2)
    \qquad \longleftrightarrow \qquad
    \Qq = \tfrac{1}{2}
    \begin{psmallmatrix}
    -m & -2r\\
    2l & m
    \end{psmallmatrix}
    \in \Sl_2(\Field)
\end{equation*}

This defines a bijective correspondence between $\Sl_2(\Field)$ and $\Qfb(\Field^2)$ through which the adjoint action of $\PGL_2(\Field)$ corresponds to the action by change of variables.
It matches the discriminants $m^2-4lr$ and sends the Lie bracket $\{\Aa,\Bb\}=\tfrac{1}{2}(\Aa\Bb-\Bb\Aa)$ of $\Sl_2(\Field)$ to the Poisson bracket of functions on $\Field^2$, under which quadratic forms are closed $\{Q_a,Q_b\} = \tfrac{1}{4}\left[(\partial_x Q_a)(\partial_y Q_b)-(\partial_x Q_b)(\partial_y Q_a)\right]
    =\{\Aa,\Bb\}$.

Consequently, all the notions defined for an element $\Qq\in \Sl_2(\Field)$ or a pair of elements $\Aa,\Bb\in \Sl_2(\Field)$, can be translated in terms of the corresponding binary quadratic forms $Q,Q_a,Q_b\in \Qfb(\Field^2)$.
For instance, after choosing a root $\sqrt{\disc(Q_a)\disc(Q_b)}=-4\sqrt{\det(\Aa\Bb)}$, we may define the cosine:
\begin{equation*}
    \cos(Q_a,Q_b)
    =\tfrac{\disc(Q_a+Q_b)-(\disc(Q_a)+\disc(Q_b))}{2\sqrt{\disc(Q_a)\disc(Q_b)}}
    =\cos(\Aa,\Bb)
\end{equation*}
and the cross-ratio of their roots $\{\alpha',\alpha\} \& \{\beta',\beta\}$, which are ordered up to simultaneous inversion:
\begin{equation*}
    \bir(Q_a,Q_b)
    =\bir(\alpha',\alpha;\beta',\beta)
    = \tfrac{(\alpha'-\alpha)(\beta'-\beta)}{(\alpha-\beta')(\beta-\alpha')}
    =\bir(\Aa,\Bb).
\end{equation*}
For a common choice of root these are related by
$\bir(Q_a,Q_b)^{-1} = \tfrac{1}{2}(1+\cos(Q_a,Q_b))$.
In particular, if $Q_a,Q_b$ have the same discriminant $\Delta$, which is to be chosen as the root $\sqrt{\disc(Q_a)\disc(Q_b)}$, then $\bir(Q_a,Q_b)^{-1} = \disc(Q_a+Q_b)/(4\Delta)$, and we may compute $\bir(Q_a,Q_b)\equiv l_al_b \bmod{\Norm_\Field \Field[\sqrt{\Delta}]}$.

The actions by change of variable of $\PGL_2(\Field)$ and $\PSL_2(\Field)$ on $\Qfb(\V)$ preserve the discriminant as an element of $\Field/(\Field^*)^2$ and $\Field$ respectively, as well as the cross ratio.
Note that the condition $\bir(Q_a,Q_b)\notin \{1,\infty\} \iff Q_a \ne \pm Q_b$ can always be achieved after conjugating $Q_b$ by $\PSL_2(\Field)$.
Consequently, Proposition \ref{Prop:Conj-a-b-K[a,b]} and Theorem \ref{Thm:Conj-a-b-SL(V)} describe the orbits of $\Qfb(\V)$ under the action of $\PGL_2(\Field)$ and $\PSL_2(\Field)$ in terms of $\disc$ and $\bir$, the latter being empty or $\Field[\sqrt{\Delta}]^\times$-torsors. 

\subsection*{The variables live in the cone.}

Conversely one may try to recover some notions defined for binary quadratic form in terms of the corresponding matrices.
This was the motivation leading to Lemma \ref{Lem:param-cone-omega}, namely to recover the values that a form $Q\in \Qfb(\V)$ takes on $\V$ in terms of the geometry of $\Qq\in \Sl(\V)$ with respect to $\Cone$.
Indeed, Lemma \ref{Lem:param-cone-omega} implies that if $\psi \colon v\in \V \mapsto p \in \Cone$ then: 
\begin{equation*}
    Q(v)= \det(v,\Qq v) = \langle \Qq,p \rangle.
\end{equation*}
So the elements $p\in \Cone$ in the cone play the role of the vector of variables $v\in \V$, whereas the other elements $\Qq\in \Sl(\V)\setminus \Cone$ are the non-degenerate binary quadratic forms $Q\in \Qfb(\V)$.

The values of $Q$ on $\V$ may thus be interpreted in terms of the distances between $\Qq^\perp$ and $\Cone$.
%
%
%
In particular, if $\disc Q_a = \Delta = \disc Q_b$ then for any $v_a,v_b\in \V$ which belong to $\Field$-bases, we have: 
\[\bir(Q_a,Q_b) \equiv Q_a(v_a)Q_b(v_b) \bmod{\Norm_\Field \Field[\sqrt{\Delta}]^\times}
\: \mathrm{whence} \:
\bir(\Aa,\Bb) \equiv \langle \Aa,p_a\rangle \langle \Bb,p_b\rangle \bmod{\Norm_\Field \Field[\sqrt{\Delta}]^\times}.
\]

\begin{Proposition} \label{Prop:ClF(D)-group}
The set $\Cl_\Field(\Delta)$ of $\PSL_2(\Field)$-orbits in $\Qfb(\V)$ with non-square discriminant $\Delta$ embeds into the group $\Field^\times /{\Norm_\Field(\Field[\sqrt{\Delta}]^\times})$ of exponent two, by sending the class of the norm $x^2-\tfrac{\Delta}{4}y^2$ of the $\Field$-extension $\Field[\sqrt{\Delta}]$ to the identity, and using the multiplication of values for composition.
\end{Proposition}


\newpage


\subsection*{Describing \texorpdfstring{$\Q$}{Q}-equivalence with Hilbert symbols}
%

In this paragraph we fix $\V=\Q^2$ and a non-square discriminant $\Delta$. Consider two forms $Q_a,Q_b$ representing variable $\Q$-classes in $\Cl_\Q(\Delta)$, and let us provide a method for computing $\bir(Q_a,Q_b)\bmod{\Norm_\Q \Q[\sqrt{\Delta}]^\times}$.
%
%

$\PP=\{-1,2\}\cup \{3,5,7\dots\}$ denotes the set of rational primes and $\Q_p$ the $p$-adic completion of $\Q$. The prime $-1$ refers (following Conway \cite{Conway_sensual-quad-form_1997}) to the place at which the completion of $\Q$ is the Archimedian field $\Q_{-1}=\R$.
%
%
For $\delta,\chi \in \Q_p^\times$ the \emph{Hilbert symbol} $\left(\delta,\chi\right)_p$ equals $1$ or $-1$ according to whether the homogenised Pell-Fermat equation $x^2-\delta y^2 = \chi z^2$ admits a solution in $\Q_p\P^2$ or not.
Thus we have $(\delta,\chi)_p = 1$ if and only if $\chi$ is the norm of an element in $\Q_p(\sqrt{\delta})$.


We define the set of prime obstructions to solving the equation $(2x)^2-\Delta y^2 = \bir(Q_a,Q_b)$ by $\PP(Q_a,Q_b)=\{p \in \PP \mid (\Delta,\bir(Q_a,Q_b))_p = -1\}$, which only depends on the $\Q$-classes of $Q_a \& Q_b$.

\begin{Theorem} The forms $Q_a$ and $Q_b$ are $\Q$-equivalent if and only if $\PP(Q_a,Q_b) = \emptyset$.
\end{Theorem}

\begin{proof}
Apply the Hasse-Minkowski theorem \cite[Chapitre IV, Th\'eor\`eme 8]{Serre_arithmetique_1970} to the ternary quadratic form $(2x)^2-\Delta y^2-\bir(Q_a,Q_b) z^2$: it represents $0$ over $\Q$ if and only if it represents $0$ over $\Q_p$ for all $p\in \PP$.
\end{proof}

The following Lemma and Remark enable us to turn the previous Theorem into a finite method for computing $\Q$-classes.

\begin{Lemma}
If $p\in \PP\setminus\{2\}$ divides $\delta$ and $\chi$ to even powers, then $(\delta,\chi)_p=1$.

In other terms $\PP(Q_a,Q_b)\setminus\{2\}$ is contained in the set of primes appearing with odd valuations in the factorisation of $\Delta$ or $\bir(Q_a,Q_b)$. In particular it is finite.
\end{Lemma}

\begin{proof}
A pedestrian method is to reduce the equation $\bmod{p}$, argue that there exists a solution by a counting procedure, and lift it to $\Q_p$ using Hensel's lemma.

Alternatively, one may use the explicit formulae \cite[Theorem III.1]{Serre_arithmetique_1970} for the Hilbert symbol at $p\ne 2$ in terms of the Legendre symbols of $\delta, \chi \in \Q_p$ at $-1$ and $p$.
\end{proof}

\begin{Remark}
Hilbert proved a global relation among the local symbols: $\prod_{p\in \PP} (\delta,\chi)_p = 1$, which is a reformulation of the quadratic reciprocity law.
Consequently $\PP(Q_a,Q_b)\setminus \{2\}$ determines $\PP(Q_a,Q_b)$.
\end{Remark}

Our final Proposition implies that $Q_a \& Q_b$ are $\Q$-equivalent if and only if $\PP(Q_0,Q_a)=\PP(Q_0,Q_b)$.
This simplifies the determination of all sets $\PP(Q_a,Q_b)$ to those involving a fixed element $Q_0$.

\begin{Proposition}
For $Q_a,Q_b,Q_c \in \Cl_\Q(\Delta)$ the set $\PP(Q_a,Q_b)$ is equal to the symmetric difference of $\PP(Q_c,Q_a)$ and $\PP(Q_b,Q_c)$.
In other terms for all $p\in \PP$ we have:
\begin{equation*}
(\Delta,\bir(Q_c,Q_a)\bir(Q_a,Q_b)\bir(Q_b,Q_c))_p = 1.
\end{equation*}
\end{Proposition}

\begin{proof}
%
According to \cite[Theorem III.2]{Serre_arithmetique_1970}, the Hilbert symbol of $\Q_p$ defines a non-degenerate symmetric bilinear form on the $\F_2$-vector space $(\Q_p^\times)/(\Q_p^\times)^2$.
The Lemma can thus be reformulated as $(\Delta,\chi_{a,b,c})_p=1$ where $\chi_{a,b,c}=\bir(Q_c,Q_a)\bir(Q_a,Q_b)\bir(Q_b,Q_c)$.
%

We must therefore compare $\chi_{a,b,c}\in \Q^\times$ with the subgroup generated by the norms of elements in $\Q(\sqrt{\Delta})^\times$.
Using the explicit formula for the cross-ratio we find that:
\begin{equation*}
    \bir(Q_a,Q_b)=\tfrac{-\Delta / (l_al_b)}{\Norm_\Delta(\alpha'-\beta)}
    \quad \mathrm{hence} \quad
    \chi_{a,b,c} =\tfrac{-\Delta^3/ (l_al_bl_c)^2}{\Norm_\Delta((\gamma'-\alpha)(\alpha'-\beta)(\beta'-\gamma))}
\end{equation*}
Consequently $(\Delta,\chi_{a,b,c})=(\Delta,-\Delta)_p = 1$ as desired.
\end{proof}

\begin{Corollary}
Denoting $Q_0$ a representative for the norm of the $\Q$-extension $\Q(\sqrt{\Delta})$, the map $Q\mapsto \PP(Q_0,Q)$ yields an isomorphism $\Cl_\Q(\Delta) = \prod_{p\in \PP} \Cl_{\Q_p}(\Delta)$.
%
\end{Corollary}


\newpage

\subsection*{Integral binary quadratic forms}

Notice that under the 1:1 correspondence $\Qfb(\Q^2) \leftrightarrow \Sl_2(\Q)$, the lattice $\Qfb(\Z^2)$ of integral binary quadratic forms gets mapped to the dual lattice $\Sl_2(\Z)^\vee$ of $\Sl_2(\Z)$ in $\Sl_2(\Q)$ with respect to the quadratic form $\det$. 
We will concentrate on the primitive points of the lattices $\Qfb(\V)$ or $\Sl(\V)^\vee$, namely those which are visible from the origin, thus not multiples of another lattice point by a non-invertible integer.

Now for any field $\Field$ of characteristic different from $2$, we may consider the extension of scalars $\Sl_2(\Z[1/2]) \to \Sl_2(\Field)$, and its restriction to $\Sl_2(\Z)^{\vee}$.
%
%
We say that $\Aa,\Bb\in \Sl_2(\Z)^{\vee}$ are $\Field$-equivalent when their images in $\Sl_2(\Field)$ belong to the same orbit under the adjoint action of $\PSL_2(\Field)$.
We may thus group the conjugacy classes of $\PSL_2(\Z)$ into $\Field$-classes, and observe how this varies with $\Field$.

When $\Field\supset \Q$ has characteristic zero, the extension of scalars $\Sl_2(\Q)\to \Sl_2(\Field)$ is injective so the $\Field$-equivalence implies the equality of discriminants.
When $\Field = \C$, this groups the integral binary quadratic forms according to their discriminant, and we find the finite class groups $\Cl(\Delta)$. 
%
When $\Field = \Q$, this defines for each discriminant $\Delta$ a partition of the class group $\Cl(\Delta)$ into $\Q$-classes.

\subsection*{Class groups and genera}

Fix a non-square discriminant $\Delta$ and consider the set $\Cl(\Delta)$ of $\PSL_2(\Z)$-equivalence classes of primitive integral binary quadratic forms with that discriminant.
%
%
This is a finite set by the classical reduction theory of binary quadratic forms (see \cite{Cox_primes-of-the-form_1997, Conway_sensual-quad-form_1997}).
%

In \cite{Gauss_disquisitiones_1807}, C.-F. Gauss endowed $\Cl(\Delta)$ with the structure of a finite abelian group 
%
%
%
which was later reformulated by Dirichlet as follows \cite{Weil_number-theory-history_1984}. One may represent two classes in $\Cl(\Delta)$ by forms $Q_a$ and $Q_b$ whose first coefficients $l_a$ and $l_b$ are coprime, and with the same middle coefficient $m$: their composition $Q_c$ of the same discriminant is determined by its first coefficient $l_c=l_al_b$ and middle coefficient $m$.
%
%
(Beware that this corresponds to the \emph{narrow} class group of ideals in $\mathcal{O}_\Delta$.)
The neutral element of $\Cl(\Delta)$ is represented by any form which takes the value $1$, for instance the principal form  $x^2+\epsilon xy+\tfrac{1}{4}(\epsilon-\Delta) y^2$ where $\epsilon \in \{0,1\}$ satisfies $\Delta=\epsilon \bmod{4}$, obtained from the $\Q$-extension $\Q(\sqrt{\Delta})$ by restricting the norm to its unique order $\mathcal{O}_\Delta$ of discriminant $\Delta$.
%

Two classes in $\Cl(\Delta)$ \emph{belong to the same genus} when for all $p\in \PP$ they are conjugate by $\SL_2(\Z_p)$, with $\Z_{-1}=\R$ (see \cite[Chapter 14]{Cassels_Rational-quadratic-forms_1978}).
%
%
%
One may consult \cite[Theorem 3.21]{Cox_primes-of-the-form_1997} for several other characterisations, such as representing the same values in $(\Z/\Delta)^\times$.
%
%
%
%
The equivalence classes for this relation form a group $\operatorname{Gen}(\Delta)$ given by the multiplication of their sets of values in $(\Z/\Delta)^\times$.
Gauss identified it with the quotient of his class group by the subgroup of squares. Moreover the kernel of the squaring map consists in the subgroup $\operatorname{Sym}(\Delta)$ of classes invariant by the Galois involution.
In other terms we have a short exact sequence of abelian groups:
    \begin{equation*}
        1\to \operatorname{Sym}(\Delta) \to \Cl(\Delta)\xrightarrow{square} \Cl(\Delta)\to \operatorname{Gen}(\Delta) \to 1.
    \end{equation*}

For two classes in $\Cl(\Delta)$ represented by $Q_a,Q_b$ we saw that $\bir(Q_a,Q_b)\equiv l_al_b \bmod{\Norm_\Q \Q[\sqrt{\Delta}]^\times}$.
%
We used this multiplication of values to define the composition law on $\Cl_\Q(\Delta)$ in proposition \ref{Prop:ClF(D)-group}.
Hence the extension of scalars $\Z\to \Q$ yields a group morphism $\Cl(\Delta) \to \Cl_\Q(\Delta)$ whose kernel consists in the $\Z$-classes which are $\Q$-equivalent to the principal form.
This kernel contains the subgroup of squares so we find a morphism of groups with exponent two $\Cl(\Delta)/\Cl(\Delta)^2 \to \Cl_\Q(\Delta)$.

A fundamental discriminant $\Delta$ is that of (the ring of integers in) a quadratic extension of $\Q$, which means that $\Delta = 1\bmod{4}$ is a square-free integer, or that $\Delta/4 \ne 1 \bmod{4}$ is a square-free integer. 
For such $\Delta$, the genus equivalence amounts to being conjugate by an element in $\PGL_2(\Q)$ (\cite[Exercise 3.17]{Cox_primes-of-the-form_1997}), which is implied by $\Q$-equivalence.
This discussion implies the following.

\begin{Proposition}
\label{Q-equiv&Genus-equiv}
In $\Cl(\Delta)$, the principal principal genus $\Cl(\Delta)^2$ is contained in the subgroup of classes which are $\Q$-equivalent to the principal form, thus genus equivalence implies $\Q$-equivalence.
%
%
If $\Delta$ is fundamental then $\Q$-equivalence implies genus equivalence (but otherwise it may not).
\end{Proposition}

\newpage

\subsection*{Example: $\Q$-equivalence and continued fractions}

Let us observe the partition of $\Cl(\Delta)=\Z/4\times \Z/2$ into $\Q$-classes for the positive fundamental discriminants $1596=4\times 3\times 7\times 19$ and $1768=8\times 13 \times 17$, whose fundamental units have norm $+1$ and $-1$.

%
%
The following tables exhibit the structure of $\Cl(\Delta) \simeq \Z/4\times \Z/2$. Each cell contains the coefficients $(l,m,r)$ of a representative $\in \Qfb(\Z^2)$ together with the period of the continued fraction expansion of its first root $\tfrac{1}{2l}(-m+\sqrt{\Delta})$ and the set of prime obstructions $\PP(Q_0,Q)$. 

\begin{table}[h]
\centering
\begin{tabular}{|c|c|c|c|}
\hline 
   $(1, -38, -38)$ 
 & $(10, -34, -11)$
 & $(25, -14, -14)$
 & $(10, -26, -23)$
 \\
   $[38, 1]$
 & $[3,1,2,3]$
 & $[1,12,1,1]$
 & $[3,3,2,1]$
 \\
   $\emptyset$
 & $\{7,19\}$ 
 & $\emptyset$
 & $\{7,19\}$
 \\
\hline
   $(2, -38, -19)$ 
 & $(5, -34, -22)$
 & $(29, -30, -6)$
 & $(5, -36, -15)$
 \\
   $[19,2]$
 & $[7,2,2,1]$
 & $[1,4,1,5]$
 & $[7,1,1,2]$
 \\
   $\{3,19\}$
 & $\{3,7\}$ 
 & $\{3,19\}$
 & $\{3,7\}$
 \\
\hline
\end{tabular}
\end{table}
\begin{table}[h]
\centering
\begin{tabular}{|c|c|c|c|}
\hline 
   $(1, 42, -1)$ 
 & $(14, 40, -3)$
 & $(-9, 34, 17)$
 & $(-3, 40, 14)$
 \\
   $[42, 42]$
 & $[2, 13, 1, 2, 13, 1]$
 & $[2, 4, 4, 2, 4, 4]$
 & $[2, 1, 13, 2, 1, 13]$
 \\
   $\emptyset$
 & $\{2,17\}$ 
 & $\emptyset$
 & $\{2,17\}$
 \\
\hline
   $(21, 40, -2)$ 
 & $(7, 40, -6)$
 & $(-13, 26, 21)$
 & $(-6, 40, 7)$
 \\
   $[1, 1, 20, 1, 1, 20]$
 & $[1, 5, 6, 1, 5, 6]$
 & $[2, 1^6, 2, 1^6]$
 & $[6, 5, 1, 6, 5, 1]$
 \\
   $\{2,13\}$
 & $\{2,17\}$ 
 & $\{2,13\}$
 & $\{2,17\}$
 \\
\hline
\end{tabular}
\end{table}

\begin{Remark}
Notice that $\Q$-equivalence does not control the period lengths of the continued fraction expansions: there exist $\Q$-equivalent forms whose roots have euclidean periods of different length.
\end{Remark}

\subsection*{Examples: $\Q$-classes modulo genera}

We may apply the methods in the previous paragraphs to determine the partition of $\Cl(\Delta)$ into $\Q$-classes (given by $\SL_2(\Q_p)$-equivalence) using the Hilbert symbols, and into genera (given by $\SL_2(\Z_p)$-equivalence) using the $\Cl(\Delta)^2$-cosets.
%
%

Denoting by $\operatorname{S}_\Q(\Delta)$ the kernel of the map $\Cl(\Delta) \to \Cl_\Q(\Delta)$, the discrepancy is measured by the dimension $c_\Q$ of the $\F_2$-vector space $\operatorname{S}_\Q(\Delta)/\Cl(\Delta)^2$.
This depends on which odd primes divide $\Delta$ to an even power, and on the $2$-adic valuation of $\Delta$.
%

In each table, the discriminant of the first row are fundamental, and the others are not.
In the second and third table, the units in the quadratic extensions have norm $+1$ and$-1$ respectively.

\begin{table}[h]
\centering
\begin{tabular}{|c|c|c|}
\hline
   $\Delta<0$
 & $\Cl(\Delta)$
 & $c_\Q$ 
 \\
\hline
   $-2^2\times 7$
 & $\Z/1$
 & $1$
 \\
\hline
   $-2^3\times 7$
 & $\Z/4$
 & $1$
 \\
\hline
    $-2^4\times 7$
 & $\Z/2$
 & $2$
 \\
 \hline
   $-2^5\times 7$
 & $\Z/4\times \Z/2$
 & $2$
 \\
\hline
   $-2^6\times 7$
 & $\Z/2\times \Z/2$
 & $4$
 \\
\hline
    $-2^7\times 7$
 & $\Z/8\times \Z/2$
 & $4$
 \\
\hline
    $-2^3\times 7^{3}$
 & $\Z/28$
 & $1$
 \\
\hline
\end{tabular}
\quad
\begin{tabular}{|c|c|c|}
\hline
   $\Delta>0$
 & $\Cl(\Delta)$
 & $c_\Q$ 
 \\
\hline
   $2^2\times 3 \times 5$
 & $\Z/2$
 & $1$
 \\
\hline
   $2^2\times 3^2 \times 5$
 & $\Z/1$
 & $1$
 \\
\hline
    $2^2\times 3 \times 5^2$
 & $\Z/2$
 & $2$
 \\
 \hline
   $2^2\times 3^3 \times 5$
 & $\Z/2$
 & $1$
 \\
\hline
   $2^2\times 3 \times 5^3$
 & $\Z/1$
 & $1$
 \\
\hline
    $2^2\times 3^3 \times 5^2$
 & $\Z/6$
 & $2$
 \\
\hline
    $2^2\times 3^2 \times 5^3$
 & $\Z/1$
 & $1$
 \\
\hline
\end{tabular}
\quad
\begin{tabular}{|c|c|c|}
\hline
   $\Delta$
 & $\Cl(\Delta)$
 & $c_\Q$ 
 \\
\hline
   $5$ \& $13$
 & $\Z/2$
 & $1$
 \\
\hline
   $5\times 13$
 & $\Z/2$
 & $1$
 \\
 \hline
    $5^2\times 13$
 & $\Z/2$
 & $2$
 \\
 \hline
 $5\times 13^2$
 & $\Z/2$
 & $2$
 \\
 \hline
 $5^3\times 13$
 & $\Z/2$
 & $1$
 \\
 \hline
 $5\times 13^3$
 & $\Z/2$
 & $1$
 \\
 \hline
 $5^3\times 13^3$
 & $\Z/2$
 & $1$
 \\
 \hline
\end{tabular}
%
\end{table}

The last example fits in the family of $\Delta = p^u\times q^v$ for distinct primes $p,q\equiv 1 \bmod{4}$ and $u,v\in \N$.
If $u=1,v=0$ then $\Delta$ is fundamental and the genera coincide with the $\Q$-classes.
For $u,v\in \N^*$ we observe that $c_\Q=1$ when both $u,v$ are odd and $c_\Q=2$ when either $u$ or $v$ is even.

\newpage

\section{Arithmetic equivalence of singular moduli and modular geodesics}
\label{Sec:modular-orbifold}

The modular group $\PSL_2(\Z)$ acts on the upper-half plane $\H\P= \{z\in \C \mid \Im(z)>0\}$ by linear fractional transformations, and the quotient is the modular orbifold $\M=\PSL_2(\Z) \backslash \H\P$.

%
Consider primitive integral binary quadratic forms $Q_a,Q_b$ with non-square discriminant $\Delta$. Fix a root $\sqrt{\Delta}$ which is to be positive if $\Delta>0$, and define the first roots of $Q_a(x,1)$ and $Q_b(x,1)$ by:
\begin{equation*}
    \alpha= \tfrac{-m_a + \sqrt{\Delta}}{2l_a} 
    \qquad \mathrm{and} \qquad
    \beta = \tfrac{-m_b + \sqrt{\Delta}}{2l_b}. 
\end{equation*} 

\subsection*{Arithmetic equivalence of singular moduli}

If $\Delta>0$, then $Q_a$ and $Q_b$ are uniquely determined by their roots $\alpha,\beta \in \H\P$.
Their $\PSL_2(\Z)$-classes correspond to points $[\alpha],[\beta]\in \M$ often called singular moduli in the study of elliptic curves.
The geodesic arc from $\alpha$ to $\beta$ in $\H\P$ has length $\lambda$ given in terms of the cross-ratio $\bir(\alpha',\alpha;\beta',\beta)$ by the formula:
\begin{equation*}
    \left(\cosh \tfrac{\lambda}{2}\right)^2
    = \frac{1 + \cosh(\lambda)}{2}
    = \frac{1}{\bir(Q_a,Q_b)}
\end{equation*}
%

\begin{Corollary}[to Theorem \ref{Thm:Conj-a-b-SL(V)}]
\label{Cor:arithmeti-equivalence-singular-moduli}
Two singular moduli $[\alpha],[\beta] \in \Field(\sqrt{\Delta})$ are $\Field$-equivalent if and only if there exists a hyperbolic geodesic arc in $\M$ from $[\alpha]$ to $[\beta]$ whose length $\lambda$ is of the form:
    \begin{equation*}
        \left(\cosh\tfrac{\lambda}{2}\right)^2 = \frac{1}{(2x)^2-\Delta y^2}
        \qquad \mathrm{for} \quad
        x,y\in \Field
    \end{equation*}
in which case all geodesic arcs from $[\alpha]$ to $[\beta]$ have this property.
\end{Corollary}

\subsection*{Arithmetic equivalence of modular geodesics}

If $\Delta < 0$, then $Q_a$ and $Q_b$ correspond to oriented geodesics $(\alpha',\alpha), (\beta',\beta)$ in $\H\P$.
Their $\PSL_2(\Z)$-classes correspond to primitive closed oriented geodesics in $\M$ called modular geodesics, whose length equals $2\sinh^{-1}(\sqrt{\Delta}/2)$.

Consider the oriented hyperbolic geodesics $(\alpha',\alpha)$ and $(\beta',\beta)$ in $\H\P$.
If they intersect, then their angle $\theta$ is given in terms of the cross-ratio $\bir(\alpha',\alpha;\beta',\beta)$ by the formula:
\begin{equation*}
    \left(\cos \tfrac{\theta}{2}\right)^2
    = \frac{1 + \cos(\theta)}{2}
    = \frac{1}{\bir(Q_a,Q_b)}
\end{equation*}
If they do not intersect, then they have a unique common perpendicular geodesic arc, which may receive compatible co-orientations from each axis or not.
When it is the case, its length $\lambda$ is given in terms of the cross-ratio $\bir(\alpha',\alpha;\beta',\beta)$ by the formula:
\begin{equation*}
    \left(\cosh \tfrac{\lambda}{2}\right)^2
    = \frac{1 + \cosh(\lambda)}{2}
    = \frac{1}{\bir(Q_a,Q_b)}
\end{equation*}


\begin{Corollary}[to Theorem \ref{Thm:Conj-a-b-SL(V)}]
\label{Cor:arithmeti-equivalence-geodesics}
Two modular geodesics of the same length $2\sinh^{-1}(\sqrt{\Delta}/2)$ are $\Field$-equivalent if and only if we have one of the following equivalent conditions:
\begin{enumerate}
    \item[$\theta$] There exists one intersection point with angle $\theta \in \,]0,\pi[$ such that:
    \begin{equation*}
        \left(\cos\tfrac{\theta}{2}\right)^2 = \frac{1}{(2x)^2-\Delta y^2} 
        \qquad \mathrm{for} \quad
        x,y\in \Field
    \end{equation*}
    in which case all intersection points have this property.
    \item[$\lambda$] There exists one co-oriented ortho-geodesic of length $\lambda$ such that:
    \begin{equation*}
        \left(\cosh\tfrac{\lambda}{2}\right)^2 = \frac{1}{(2x)^2-\Delta y^2}
        \qquad \mathrm{for} \quad
        x,y\in \Field
    \end{equation*}
    in which case all co-oriented ortho-geodesics have this property.
\end{enumerate}
\end{Corollary}

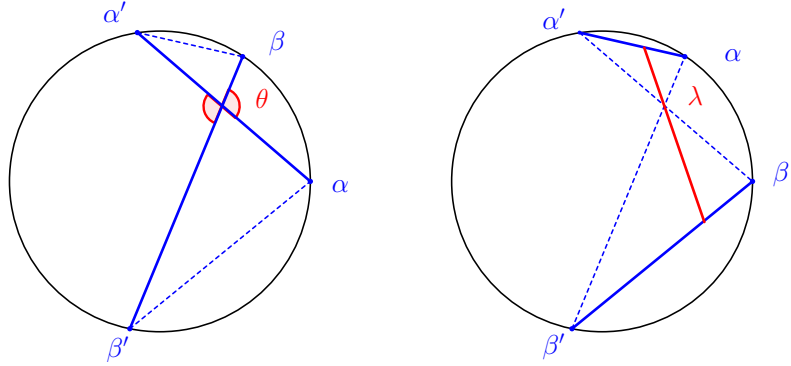
\begin{figure}[h]
    \centering
    \scalebox{.5}{
\begin{tikzpicture}[line cap=round,line join=round,>=triangle 45,x=4.0cm,y=4.0cm]
\clip(-1.2,-1.2) rectangle (1.3,1.3);
\draw [shift={(0.41,0.5)},line width=1.8pt,color=red,fill=red,fill opacity=0.1] (0,0) -- (-40.55:0.12) arc (-40.55:67.61:0.12) -- cycle;
\draw [shift={(0.41,0.5)},line width=1.8pt,color=red,fill=red,fill opacity=0.1] (0,0) -- (139.45:0.12) arc (139.45:247.61:0.12) -- cycle;
\draw [line width=1.2pt] (0,0) circle (4cm);
\draw [line width=2pt,color=blue] (-0.15,0.99)-- (1,0);
\draw [line width=2pt,color=blue] (-0.2,-0.98)-- (0.55,0.83);
\draw [line width=1.2pt,dash pattern=on 3pt off 3pt,color=blue] (-0.15,0.99)-- (0.55,0.83);
\draw [line width=1.2pt,dash pattern=on 3pt off 3pt,color=blue] (-0.2,-0.98)-- (1,0);
\begin{scriptsize}
\fill [color=blue] (-0.15,0.99) circle (2.0pt);
\draw[color=blue] (-0.3,1.13) node {\Huge$\alpha'$};
\fill [color=blue] (1,0) circle (2.0pt);
\draw[color=blue] (1.2,-0.03) node {\Huge$\alpha$};
\fill [color=blue] (-0.2,-0.98) circle (2.0pt);
\draw[color=blue] (-0.28,-1.12) node {\Huge$\beta'$};
\fill [color=blue] (0.55,0.83) circle (2.0pt);
\draw[color=blue] (0.78,0.92) node {\Huge$\beta$};
\draw[color=red] (0.68,0.55) node {\Huge$\theta$};
\end{scriptsize}
\end{tikzpicture}
}
    \qquad
    \scalebox{.5}{
\begin{tikzpicture}[line cap=round,line join=round,>=triangle 45,x=4.0cm,y=4.0cm]
\clip(-1.2,-1.2) rectangle (1.3,1.3);
\draw [line width=1.2pt] (0,0) circle (4cm);
\draw [line width=1.2pt,dash pattern=on 3pt off 3pt,color=blue] (-0.15,0.99)-- (1,0);
\draw [line width=1.2pt,dash pattern=on 3pt off 3pt,color=blue] (-0.2,-0.98)-- (0.55,0.83);
\draw [line width=2pt,color=blue] (-0.15,0.99)-- (0.55,0.83);
\draw [line width=2pt,color=blue] (-0.2,-0.98)-- (1,0);
\draw [line width=2pt,color=red] (0.68,-0.27)-- (0.28,0.89);
\begin{scriptsize}
\fill [color=blue] (-0.15,0.99) circle (1.5pt);
\draw[color=blue] (-0.32,1.08) node {\Huge$\alpha'$};
\fill [color=blue] (1,0) circle (2.0pt);
\draw[color=blue] (1.19,0.05) node {\Huge$\beta$};
\fill [color=blue] (-0.2,-0.98) circle (2.0pt);
\draw[color=blue] (-0.32,-1.1) node {\Huge$\beta'$};
\fill [color=blue] (0.55,0.83) circle (2.0pt);
\draw[color=blue] (0.87,0.84) node {\Huge$\alpha$};
\draw[color=red] (0.62,0.57) node {\Huge$\lambda$};
\end{scriptsize}
\end{tikzpicture}
}
    \caption*{Cross-ratios and cosines in the real case.}
\end{figure}

\begin{figure}[h]
    \centering
    \scalebox{.3}{
\begin{tikzpicture}[line cap=round,line join=round,>=triangle 45,x=1.0cm,y=1.0cm]
\clip(-4.,-4.) rectangle (4.,4.);
\draw [shift={(0.,0.)},line width=2.pt,color=red,fill=red,fill opacity=0.10000000149011612] (0,0) -- (56.309932474020215:2.3302255765350153) arc (56.309932474020215:123.6900675259798:2.3302255765350153) -- cycle;
\draw [->,line width=2.pt,color=blue] (-2.,-3.) -- (2.,3.);
\draw [->,line width=2.pt,color=blue] (2.,-3.) -- (-2.,3.);
\begin{scriptsize}
\draw[color=blue] (-2.413912953651292,-2.8340125566941947) node {\Huge$\alpha'$};
\draw[color=blue] (2.9223036166138927,3.084760407704745) node {\Huge$\alpha$};
\draw[color=blue] (2.9223036166138927,-2.880617068224895) node {\Huge$\beta'$};
\draw[color=blue] (-2.4838197209473423,3.084760407704745) node {\Huge$\beta$};
\draw [fill=red] (0.,0.) circle (2.0pt);
\draw[color=red] (0.26584645936397566,3.014853640408694) node {\Huge$\theta$};
\end{scriptsize}
\end{tikzpicture}
}
    \qquad 
    \scalebox{.3}{
\begin{tikzpicture}[line cap=round,line join=round,>=triangle 45,x=1.0cm,y=1.0cm]
\clip(-4.,-4.) rectangle (4.,4.);
\draw[line width=2.pt,color=violet,fill=violet,fill opacity=0.10000000149011612] (2.,0.3772399719250914) -- (1.6227600280749086,0.37723997192509146) -- (1.6227600280749086,0.) -- (2.,0.) -- cycle; 
\draw[line width=2.pt,color=violet,fill=violet,fill opacity=0.10000000149011612] (-1.6227600280749086,0.) -- (-1.6227600280749086,0.37723997192509134) -- (-2.,0.3772399719250914) -- (-2.,0.) -- cycle; 
\draw [->,line width=2.pt,color=blue] (-2.,-3.) -- (-2.,3.);
\draw [->,line width=2.pt,color=blue] (2.,-3.) -- (2.,3.);
\draw [line width=2.pt,color=red] (-2.,0.)-- (2.,0.);
\draw [->,line width=2.pt,color=red] (-1.,0.) -- (-1.,1.);
\draw [->,line width=2.pt,color=red] (1.,0.) -- (1.,1.);
\begin{scriptsize}
\draw[color=blue] (-2.6,-2.7) node {\Huge$\alpha'$};
\draw[color=blue] (-2.6,2.7) node {\Huge$\alpha$};
\draw[color=blue] (2.6,2.7) node {\Huge$\beta'$};
\draw[color=blue] (2.6,-2.7) node {\Huge$\beta$};
\end{scriptsize}
\end{tikzpicture}
\begin{tikzpicture}[line cap=round,line join=round,>=triangle 45,x=1.0cm,y=1.0cm]
\clip(-4.,-4.) rectangle (4.,4.);
\draw[line width=2.pt,color=violet,fill=violet,fill opacity=0.10000000149011612] (1.6227600280749086,0.) -- (1.6227600280749086,-0.3772399719250923) -- (2.,-0.3772399719250923) -- (2.,0.) -- cycle; 
\draw[line width=2.pt,color=violet,fill=violet,fill opacity=0.10000000149011612] (-1.6227600280749086,0.) -- (-1.6227600280749086,0.37723997192509134) -- (-2.,0.3772399719250914) -- (-2.,0.) -- cycle; 
\draw [->,line width=2.pt,color=blue] (-2.,-3.) -- (-2.,3.);
\draw [->,line width=2.pt,color=blue] (2.,3.) -- (2.,-3.);
\draw [line width=2.pt,color=red] (-2.,0.)-- (2.,0.);
\draw [->,line width=2.pt,color=red] (-1.,0.) -- (-1.,1.);
\draw [->,line width=2.pt,color=red] (1.,0.) -- (1.,-1.);
\begin{scriptsize}
\draw[color=blue] (-2.6,-2.7) node {\Huge$\alpha'$};
\draw[color=blue] (-2.6,2.7) node {\Huge$\alpha$};
\draw[color=blue] (2.6,2.7) node {\Huge$\beta'$};
\draw[color=blue] (2.6,-2.7) node {\Huge$\beta$};
\end{scriptsize}
\end{tikzpicture}
}
    \caption*{Angle well defined in $\,]0,\pi[$. Ortho-geodesics well and badly co-oriented.}
\end{figure}
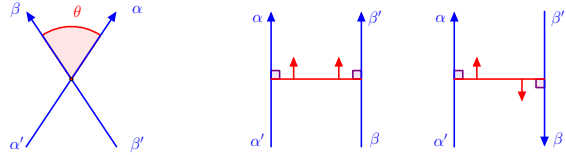

\begin{Remark}
    Since Theorem \ref{Thm:Conj-a-b-SL(V)} holds over any integral ring in which $2$ is invertible, one may similarly characterise genus equivalence (for all discriminants $\Delta$) replacing $\Q$ with $\Z_p$ for all $p\in \PP$. 
\end{Remark}


\subsection*{Linking numbers of modular knots}

Let us briefly relate the arithmetic-geometric intersections of modular geodesics to the topological linking numbers of modular knots,  referring to \cite{CLS_Link-PSL2K_2022} as well as \cite[Chapter 0 and Chapter 5]{CLS_phdthesis_2022} for the details and much more about this.

The unit tangent bundle of the modular orbifold $\M = \PSL_2(\Z)\backslash \H\P$ can be identified with the manifold $\U = \PSL_2(\Z)\backslash \PSL_2(\R)$, homeomorphic to the complement of a trefoil knot in the sphere.
The primitive closed geodesics of $\M$ lift in $\U$ to the primitive periodic orbits for the geodesic flow: one may ask about the linking numbers between these \emph{modular knots}. 
They correspond to the $\PSL_2(\Z)$-classes of primitive $A\in \PSL_2(\Z)$, or of primitive $\Aa \in \Sl_2(\Z)^\vee$, with positive discriminant.

\begin{figure}[h]
    \centering
    \includegraphics[width=0.36\textwidth]{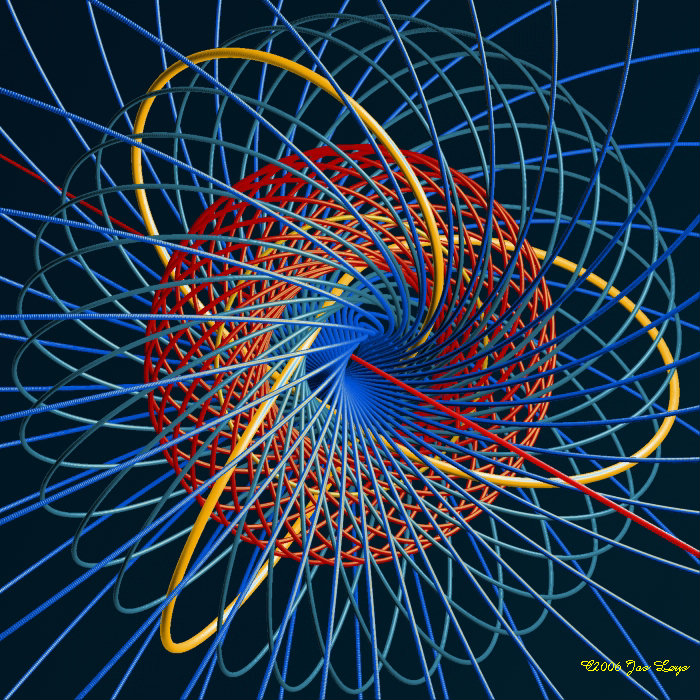}
    \includegraphics[width=0.48\textwidth]{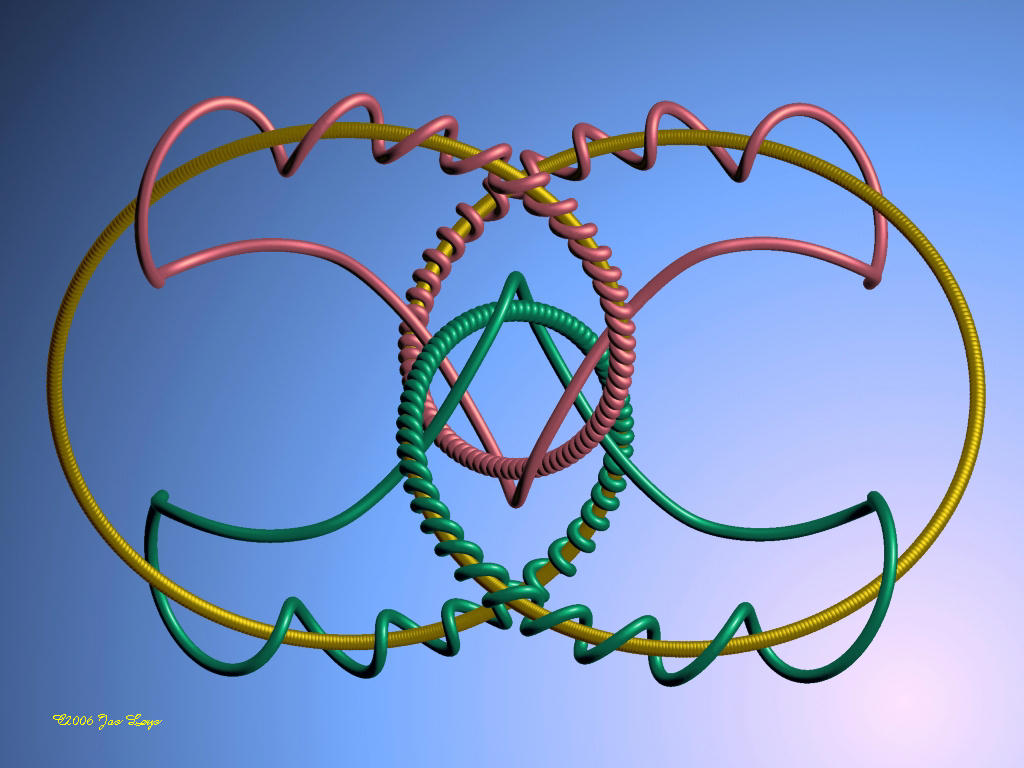}
    \caption*{
    The Seifert fibration $\U\to \M$ and two modular knots, from the \href{http://www.josleys.com/articles/ams_article/Lorenz3.htm}{online article} \cite{GhyLey_Lorenz-Modular-visual_2016} which proposes an animated introduction to the topology and dynamics of $\U$.
    }
\end{figure}

Let us introduce, for any pair of modular geodesics $\gamma_A,\gamma_B$, the following sums of the arithmetic-geometric quantities encountered in Corollary \ref{Cor:arithmeti-equivalence-geodesics} over their oriented intersection angles $\theta \in \,]0,\pi[$:
\begin{equation*}
    \Link_q(A,B)
    = \tfrac{1}{2} \sum \left(\cos \tfrac{\theta}{2}\right)^2 
    \qquad \mathrm{and} \qquad
    \Cos_q(A,B)
    = \tfrac{1}{2} \sum \left(\cos \theta\right)
\end{equation*}
and study their variations as we deform the metric on $\M$ by opening the cusp.

The complete hyperbolic metrics on the orbifold $\M$ correspond to the faithful and discrete representations $\rho\colon \PSL_2(\Z) \to \PSL_2(\R)$ up to conjugacy.
Since $\PSL_2(\Z)$ is the free amalgam of its cyclic subgroups of order $2$ and $3$ generated by $S = \begin{psmallmatrix}0 & -1 \\ 1 & 0 \end{psmallmatrix}$ and $T = \begin{psmallmatrix} 1 & -1 \\ 1 & 0 \end{psmallmatrix}$, they form a $1$-dimensional real algebraic set parametrized by $q\in \R_+^*$, fixing $S$ and conjugating $T$ by $\exp\left(-\tfrac{1}{2}\log(q)K\right)$:
%
\begin{equation*}
    S_q =
    \begin{pmatrix}
    0 & -1 \\
    1 & 0
    \end{pmatrix}
    \quad
    T_q = 
    \begin{pmatrix}
    1 & -q \\
    q^{-1} & 0
    \end{pmatrix}
    \qquad 
    L_q = 
    \begin{pmatrix}
    q & 0 \\
    1 & q^{-1}
    \end{pmatrix}
    \quad
    R_q = 
    \begin{pmatrix}
    q & 1 \\
    0 & q^{-1}
    \end{pmatrix}
\end{equation*}
%
The geometric-algebra of $\PSL_2(\Field)$ applies in particular to the image of $\rho_q \colon \PSL_2(\Z) \to \PSL_2(\R)$. For example, if the axes of $A_q, B_q \in \SL_2(\R)$ intersect, the cosine of their oriented angle is given by:
\begin{equation*}
    \cos(A_q,B_q)
    = \sign(\Tr(A_q) \Tr(B_q))\frac{\Tr(A_qB_q)-\Tr(A_qB_q^{-1})}{\sqrt{\disc(A_q)\disc(B_q)}}.
\end{equation*}

The primitive hyperbolic conjugacy classes of $\PSL_2(\Z)$ still index the hyperbolic geodesics in the quotient $\M_q=\rho_q(\PSL_2(\Z))\backslash\H\P$ which do not surround the cusp.
We may thus define the analogous sums $\Link_q(A,B)$ and $\Cos_q(A,B)$ over the intersection angles $\theta_q \in \,]0,\pi[$ between the $q$-modular geodesics $\gamma_{A_q},\gamma_{B_q} \subset \M_q$ of the $\tfrac{1}{2}\left(\cos \tfrac{1}{2}\theta_q\right)^{2}$ and $\left(\cos \theta_q\right)$.

As $q\to \infty$, the hyperbolic orbifold $\M_q$ has a convex core which retracts onto a thin neighbourhood of the long geodesic arc connecting its conical singularities, whose preimage in the universal cover $\H\P$ is a trivalent tree. In the limit we recover the action of $\PSL_2(\Z)$ on its Bruhat-Tits building, the infinite planar trivalent tree $\TT$, and by studying its combinatorics \cite{CLS_Link-PSL2K_2022} proves the following.

\begin{Theorem}[Linking and intersection from boundary evaluations]
\label{Thm:Bir(A,B)-->lk(A,B)}
For primitive hyperbolic $A,B\in \PSL_2(\Z)$, the limits of the function $\Link_q(A,B)$ and $\Cos_q(A,B)$ at the boundary point of the $\PSL_2(\R)$-character variety of $\PSL_2(\Z)$ recover their linking and intersection numbers:
\begin{align*}
    &\Link_q(A,B) \xrightarrow[q\to \infty]{} \lk(A,B)
    \\
    &\Cos_q(A,B) \xrightarrow[q\to \infty]{} 
    \lk(A,B)-\lk(A^{-1},B)
    =\lk(A,B)-\tfrac{1}{4}I(A, B)
\end{align*}
\end{Theorem}

Hence the functions $\Link_q \& \Cos_q$ interpolate between the geometry at $q=1$ of the arithmetic group $\PSL_2(\Z) \subset \PSL_2(\R)$ and the topology at $q=+\infty$ of the combinatorial action $\PSL_2(\Z) \to \Aut(\TT)$.

\begin{Remark}
This discussion naturally carries over to the field of functions $\Field=\Q(q)$ on the character variety of $\PSL_2(\Z)$, or its universal quadratic closure.

The series $L_q(A,B)$ is thus reminiscent of the special value at $s=2$ of a restricted zeta function: its terms $\tfrac{1}{2}\left(1+\cos(A_q,B_q)\right)=1/\bir(A_q, B_q)$ are the inverse norms of certain principal ideals for the quadratic extension of $\Q(q)$ generated by $\sqrt{\disc(A_q)\disc(B_q)}$.

\end{Remark}

Let us display the graphs of $\textcolor{blue}{q\mapsto2\Link_q(A,B)}$ and $\textcolor{red}{q\mapsto2\Link_q(A,B^{-1})}$ along with their average $\textcolor{black!50!green}{\tfrac{1}{2}I(A,B)}$ for some pairs $A,B\in \PSL_2(\N)$. The legend $A=[a_0,a_1,\dots]$ means $A=R^{a_0}L^{n_1}\dots$ has attractive fixed point $\alpha \in \R\P^1$ with continued fraction $\alpha=a_0+\tfrac{1}{a_1+\dots}$.
%


\begin{figure}[h]
    \centering
    \includegraphics[width=0.36\textwidth]{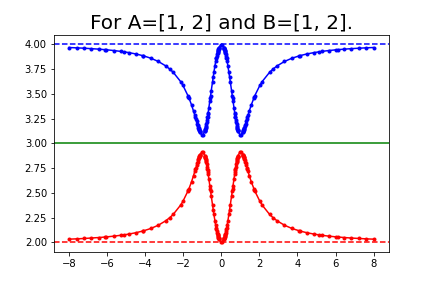}
    \hspace{-0.9cm}
    \includegraphics[width=0.36\textwidth]{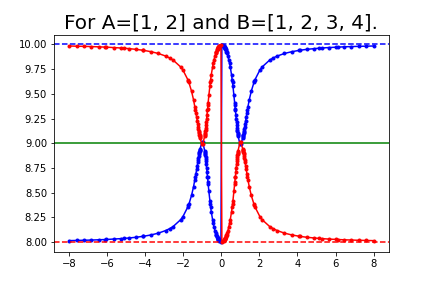}
    \hspace{-0.9cm}
    \includegraphics[width=0.36\textwidth]{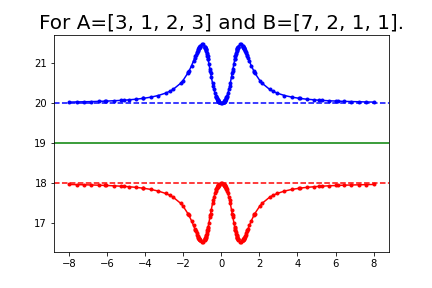}
    \vspace{-0.4cm}
    \caption*{\textcolor{blue}{$L_q(A,B)$} interpolates between the arithmetic at $1$ and the topology at $+\infty$.}
\end{figure}

\begin{Remark}
Variation on Theorem \ref{Thm:Bir(A,B)-->lk(A,B)} can be obtained by applying by any continuous function on the interval $[0,1]$ to the terms of $\Link_q(A,B)$.
A particularly interesting example is given by
\begin{equation*}
    \Link_q(A,B)
    = \tfrac{1}{2} \sum \LL \left(\cos \tfrac{\theta}{2}\right)^2 
    = \tfrac{1}{2} \sum \LL\left( \bir(\alpha',\alpha;\beta',\beta)\right)
\end{equation*}
where $\LL(z)=\sum_{n=1}^{\infty} \frac{z^n}{n^2}+\tfrac{1}{2} \log(\lvert z \rvert)\log(1-z)$ is Rogers' normalisation of the dilogarithm, yielding:
\begin{align*}
    &\LL_q(A,B) \xrightarrow[q\to \infty]{} \tfrac{\pi}{6}\lk(A,B)
\end{align*}

\end{Remark}

\bibliographystyle{alpha} 
\bibliography{biblio-these.bib}

\begin{thebibliography}{Sim22b}

\bibitem[Arn05]{Arnold_Jacobi-Lie_2005}
V.~Arnold.
\newblock Lobachevsky triangle altitudes theorem as the {J}acobi identity in
  the {L}ie algebra of quadratic forms on symplectic plane.
\newblock {\em J. Geom. Phys.}, 53(4):421--427, 2005.

\bibitem[Car92]{Cartan_Lecons-3-reprint_1992}
\'{E}lie Cartan.
\newblock {\em Le\c{c}ons sur la g\'{e}om\'{e}trie projective complexe. {L}a
  th\'{e}orie des groupes finis et continus et la g\'{e}om\'{e}trie
  diff\'{e}rentielle trait\'{e}es par la m\'{e}thode du rep\`ere mobile.
  {L}e\c{c}ons sur la th\'{e}orie des espaces \`a connexion projective}.
\newblock Les Grands Classiques Gauthier-Villars. [Gauthier-Villars Great
  Classics]. \'{E}ditions Jacques Gabay, Sceaux, 1992.
\newblock Reprint of the editions of 1931, 1937 and 1937.

\bibitem[Cas78]{Cassels_Rational-quadratic-forms_1978}
J.~W.~S. Cassels.
\newblock {\em Rational quadratic forms}, volume~13 of {\em London Mathematical
  Society Monographs}.
\newblock Academic Press, Inc. [Harcourt Brace Jovanovich, Publishers],
  London-New York, 1978.

\bibitem[CF97]{Conway_sensual-quad-form_1997}
John Conway and Francis Fung.
\newblock {\em The sensual (quadratic) form}.
\newblock MAA, 1997.

\bibitem[Cox97]{Cox_primes-of-the-form_1997}
David Cox.
\newblock {\em Primes of the form $x^2 + ny^2$}.
\newblock Wiley-Interscience, 1997.

\bibitem[Die71]{Dieudonne_geometrie-groupes_1971}
Jean Dieudonn\'e.
\newblock {\em La g\'eometrie des groupes classiques}.
\newblock Springer-Verlag, 3 edition, 1971.

\bibitem[Gau07]{Gauss_disquisitiones_1807}
Carl~Friedrich Gauss.
\newblock {\em Recherches Arithm\'etiques}.
\newblock Courcier, 1807.

\bibitem[GL16]{GhyLey_Lorenz-Modular-visual_2016}
\'{E}tienne Ghys and Jos Leys.
\newblock Lorenz and modular flows: a visual introduction.
\newblock http://www.ams.org/publicoutreach/feature-column/fcarc-lorenz, 2016.

\bibitem[Iva11]{Ivanov_Arnold-Jacobi-Lie_2011}
Nikolai Ivanov.
\newblock Arnol'd, the {J}acobi identity, and orthocenters.
\newblock {\em Amer. Math. Monthly}, 118(1):41--65, 2011.

\bibitem[Pen87]{Penner_deco-teich-space-punct-surf_1987}
Robert Penner.
\newblock The decorated {T}eichm\"{u}ller space of punctured surfaces.
\newblock {\em Comm. Math. Phys.}, 113(2):299--339, 1987.

\bibitem[Pen96]{Penner_geometry-gauss_1996}
Robert Penner.
\newblock The geometry of the {G}auss product.
\newblock {\em J. Math. Sci.}, 81(3):2700--2718, 1996.

\bibitem[Ser70]{Serre_arithmetique_1970}
Jean-Pierre Serre.
\newblock {\em Cours d'arithm\'etique}.
\newblock PUF, 1970.

\bibitem[Sim22a]{CLS_phdthesis_2022}
Christopher-Lloyd Simon.
\newblock {\em {Arithmetic and Topology of Modular knots}}.
\newblock Th{\`e}se, {Universit{\'e} de Lille}, June 2022.
\newblock
  \href{https://tel.archives-ouvertes.fr/tel-03755147/file/AriTopoModuKnots_14-07-2022.pdf}{PDF
  on HAL}.

\bibitem[Sim22b]{CLS_Link-PSL2K_2022}
Christopher-Lloyd Simon.
\newblock Linking numbers of modular knots, 2022.
\newblock Submitted for publication,
  \href{https://arxiv.org/abs/2211.05957}{arxiv version}.

\bibitem[Wei84]{Weil_number-theory-history_1984}
Andr\'{e} Weil.
\newblock {\em Number theory: {A}n approach through history From {H}ammurapi to
  {L}egendre}.
\newblock Birkh\"{a}user, 1984.

\end{thebibliography}

\end{document}